\documentclass{imsart}
\RequirePackage{amsthm,amsmath,amsfonts,amssymb}
\RequirePackage[numbers]{natbib}
\RequirePackage[colorlinks,citecolor=blue,urlcolor=blue]{hyperref}
\RequirePackage{graphicx,ulem}
\usepackage{booktabs}
\usepackage[noend]{algpseudocode}
\usepackage{algorithmicx,algorithm}
\usepackage{fancyhdr}
\startlocaldefs

\newtheorem{theorem}{Theorem}[section]

\theoremstyle{remark}

\endlocaldefs

\begin{document}

\newcommand{\bbA}{{\bf A}}
\newcommand{\cbbA}{{\check \bbA}}
\newcommand{\mbbA}{{\mathcal A}}
\newcommand{\bbbA}{{\bar \mbbA}}
\newcommand{\hbbA}{{\widehat \bbA}}
\newcommand{\tbbA}{{\tilde \bbA}}
\newcommand{\bba}{{\bf a}}
\newcommand{\bbB}{{\bf B}}
\newcommand{\bbC}{{\bf C}}
\newcommand{\bbc}{{\bf c}}
\newcommand{\bbD}{{\bf D}}
\newcommand{\bbd}{{\bf d}}
\newcommand{\bbe}{{\bf e}}
\newcommand{\bbE}{{\bf E}}
\newcommand{\rE}{{\rm E}}
\newcommand{\bbf}{{\bf f}}
\newcommand{\bbF}{{\bf F}}
\newcommand{\bbP}{{\bf P}}
\newcommand{\bbg}{{\bf g}}
\newcommand{\bbG}{{\bf G}}
\newcommand{\bbK}{{\bf K}}
\newcommand{\bbH}{{\bf H}}
\newcommand{\bbh}{{\bf h}}
\newcommand{\bbw}{{\bf w}}
\newcommand{\bbI}{{\bf I}}
\newcommand{\bbi}{{\bf i}}
\newcommand{\bbj}{{\bf j}}
\newcommand{\bbJ}{{\bf J}}
\newcommand{\bbk}{{\bf k}}
\newcommand{\bbl}{{\bf 1}}
\newcommand{\bbM}{{\bf M}}
\newcommand{\bbm}{{\bf m}}
\newcommand{\bbN}{{\bf N}}
\newcommand{\bbn}{{\bf n}}
\newcommand{\bbQ}{{\bf Q}}
\newcommand{\bbq}{{\bf q}}
\newcommand{\bbO}{{\bf O}}
\newcommand{\bbR}{{\bf R}}
\newcommand{\bbr}{{\bf r}}
\newcommand{\bbs}{{\bf s}}
\newcommand{\cbbs}{{\check \bbs}}
\newcommand{\hbbs}{{\hat \bbs}}
\newcommand{\tbbs}{{\tilde \bbs}}
\newcommand{\bbS}{{\bf S}}
\newcommand{\cbbS}{{\check \bbS}}
\newcommand{\wB}{{\widehat B}}
\newcommand{\tbbS}{\widetilde{\bf S}}
\newcommand{\hbbS}{\widehat{\bf S}}
\newcommand{\obbS}{\overline{\bf S}}
\newcommand{\bbt}{{\bf t}}
\newcommand{\bbT}{{\bf T}}
\newcommand{\hbbT}{\widetilde{\bf T}}
\newcommand{\tbbT}{\widetilde{\bf T}}
\newcommand{\obT}{{\overline{\bf T}}}
\newcommand{\bbU}{{\bf U}}
\newcommand{\bbu}{{\bf u}}
\newcommand{\bbV}{{\bf V}}
\newcommand{\bbv}{{\bf v}}
\newcommand{\tbbv}{\widetilde{\bf v}}
\newcommand{\bbW}{{\bf W}}
\newcommand{\tbbW}{\widetilde{\bf W}}
\newcommand{\hbbW}{\widehat{\bf W}}
\newcommand{\bbX}{{\bf X}}
\newcommand{\tbby}{\tilde {\bf y}}
\newcommand{\tbbX}{\widetilde {\bf X}}
\newcommand{\hbbX}{\widehat {\bf X}}
\newcommand{\bbx}{{\bf x}}
\newcommand{\obbx}{{\overline{\bf x}}}
\newcommand{\tbbx}{{\widetilde{\bf x}}}
\newcommand{\hbbx}{{\widehat{\bf x}}}
\newcommand{\obby}{{\overline{\bf y}}}
\newcommand{\hbbY}{{\widehat{\bf Y}}}
\newcommand{\tbbY}{{\widetilde{\bf Y}}}
\newcommand{\bbY}{{\bf Y}}
\newcommand{\bby}{{\bf y}}
\newcommand{\bbZ}{{\bf Z}}
\newcommand{\bbz}{{\bf z}}
\newcommand{\bbb}{{\bf b}}
\newcommand{\cD}{{\cal D}}
\newcommand{\bbL}{{\bf L}}
\newcommand{\bxi} {\boldsymbol  \xi}
\newcommand{\bbeta} {\boldsymbol  \eta}
\newcommand{\utm}{\underline{ \tilde m}}
\newcommand{\um}{\underline{m}}

\newcommand{\rdd}{\textcolor{red}}
\newcommand{\bll}{\textcolor{blue}}

\newcommand{\md}{\mbox{d}}
\newcommand{\non}{\nonumber\\}
\newcommand{\tr}{{\rm tr}}
\newcommand{\E}{{\mathbb{E}}}
\newcommand{\rP}{{\mathbb{P}}}
\newcommand{\bqa}{\begin{eqnarray}}
\newcommand{\eqa}{\end{eqnarray}}

\newcommand{\bqn}{\begin{eqnarray*}}
\newcommand{\eqn}{\end{eqnarray*}}

\newtheorem{thm}{Theorem}[section]
\newtheorem{corollary}[thm]{Corollary}
\newtheorem{defin}[thm]{Definition}
\newtheorem{prop}[thm]{Proposition}
\newtheorem{remark}[thm]{Remark}
\newtheorem{assumption}[thm]{Assumption}

\begin{frontmatter}
\title{Spiked eigenvalues of noncentral Fisher matrix with applications\thanks{The first two authors contributed equally to this work. For correspondence, please contact Zhidong Bai and Jiang Hu.}}
\runtitle{Spiked eigenvalues of noncentral Fisher matrix with applications}

\begin{aug}
	\author[A]{\fnms{Xiaozhuo} \snm{Zhang}\ead[label=e1]{zhangxz722@nenu.edu.cn}},
	\author[B]{\fnms{Zhiqiang} \snm{Hou}\ead[label=e2]{houzq399@nenu.edu.cn}},
	\author[A]{\fnms{Zhidong} \snm{Bai}\ead[label=e3]{baizd@nenu.edu.cn}}
	\and
	\author[A]{\fnms{Jiang} \snm{Hu}\ead[label=e4]{huj156@nenu.edu.cn}}
	
	\address[A]{School of Mathematics and Statistics, Northeast Normal University,
		\printead{e1,e3,e4}}
	
	\address[B]{School of Statistics, Shandong University of Finance and Economics,
		\printead{e2}}
\end{aug}
\runauthor{X.Z. Zhang, Z.Q. Hou, Z.D. Bai, J. Hu}

\begin{abstract}
In this paper, we investigate the asymptotic behavior of spiked eigenvalues of the noncentral Fisher matrix defined by ${\mathbf F}_p={\mathbf C}_n(\mathbf S_N)^{-1}$,	where ${\mathbf C}_n$ is a noncentral sample covariance matrix defined by $(\mathbf \Xi+\mathbf X)(\mathbf \Xi+\mathbf X)^*/n$ and $\mathbf S_N={\mathbf Y}{\mathbf Y}^*/N$. The matrices $\mathbf X$ and $\mathbf Y$ are two independent {Gaussian} arrays, with respective $p\times n$ and $p\times N$ and the Gaussian entries of them are \textit {independent and identically distributed} (i.i.d.) with mean $0$ and variance $1$. When $p$, $n$, and $N$ grow to infinity proportionally, we establish a phase transition of the spiked eigenvalues of $\mathbf F_p$. Furthermore, we derive the \textit{central limiting theorem} (CLT) for the spiked eigenvalues of $\mathbf F_p$. As an accessory to the proof of the above results, the fluctuations of the spiked eigenvalues of ${\mathbf C}_n$ are studied, which should have its own interests. Besides,  we develop the limits and CLT for the sample canonical correlation coefficients by the results of the spiked noncentral Fisher matrix and  give three consistent estimators, including the population spiked eigenvalues and the population canonical correlation coefficients.
\end{abstract}

\end{frontmatter}

\pagestyle{myheadings}
\chead{\runtitle{Spiked eigenvalues of noncentral Fisher matrix with applications}} 
\chead{\runauthor{X.Z. Zhang, Z.Q. Hou, Z.D. Bai, J. Hu}}

\section{Introduction}
Fisher matrix is one of the most classical and important tools in multivariate statistic analysis (for details see \cite{Anderson2003}, \cite{Ma1979}, and \cite{Mu1982}). \cite{J1964} provided a remarked five-way classification of the distribution theory and introduced some representative applications, such as signal detection in noise and testing equality of group means under unknown covariance matrix and so on. Among these applications, some statistics can be transformed into a Fisher matrix while others can be studied by a noncentral Fisher matrix. So it is natural to study the spectral properties of the Fisher matrix and noncentral Fisher matrix. 

There have been many works focusing on the Fisher matrix, \cite{W1980} derived the \textit{limiting spectral distribution} (LSD) of Fisher matrix, which is the celebrated \textit{Wachter distribution}.  \cite{HPZ2016} proved the largest eigenvalue of Fisher matrix follows \textit{Tracy-Widom} (T-W) law, see \cite{TW1996}. \cite{Zheng2012} was devoted to the CLT for \textit{linear spectral statistics} (LSS) of the Fisher matrix. \cite{ZBY2017} studied the LSD and {CLT of} LSS of the so-called general Fisher matrix. {In fact}, these above works all focus on the central Fisher matrix. Before introducing the concept of the noncentral Fisher {matrix}, it is necessary to know the large dimensional information-plus-noise-type matrix,
\begin{align}\label{noncentralS}
	{\mathbf C}_n=\frac{1}{n}(\mathbf \Xi+\mathbf X)(\mathbf \Xi+\mathbf X)^*
\end{align}
where $\mathbf X$ is a $p\times n$ matrix containing i.i.d. entries with mean $0$ and variance $1$.  ${\mathbf \Xi}{\mathbf \Xi}^*/n$ is {a} deterministic  and it is assumed to have a LSD. {Here and subsequently,  $*$ denotes conjugate transpose, and  $T$ stands for transpose on real matrix and vector.} In \cite{DS2007a,DS2007b,C2014,BS2012}, a lot of spectral properties of ${\mathbf C}_n$ have been researched. Actually, the matrix ${\mathbf C}_n$ is a noncentral sample covariance matrix and $\boldsymbol{\Xi}\boldsymbol{\Xi}^{\ast}/n$ is called the noncentral parameter matrix, whose eigenvalues are arranged as a descending order:
\begin{align}\label{Zhuaneig3}
	l_1^{\boldsymbol{\Xi}}\geq l_2^{\boldsymbol{\Xi}} \geq \cdots \geq l_p^{\boldsymbol{\Xi}}.
\end{align}  However, many problems, such as {signal detection in noise} and testing equality of group means under unknown covariance matrix always involve {the} noncentral Fisher matrix, which is constructed based on the matrix ${\mathbf C}_n$,
\begin{align}\label{noncentralF}
	{\mathbf F}_p={\mathbf C}_n{\mathbf S}_N^{-1},
\end{align}
where  ${\mathbf S}_N={\mathbf Y}{\mathbf Y}^*/N$ and $\mathbf Y$  is independent of $\mathbf X$. The entries  $\{\mathbf Y_{ij}, 1\leq i\leq p, 1\leq j \leq N\}$  are i.i.d. with mean $0$ and variance $1$. To the best of our knowledge, there are only a handful of works devoted to the noncentral Fisher matrices. Under the Gaussian assumption, \cite{JN2017} developed an approximation to the distribution of the largest eigenvalue of the noncentral Fisher matrix and \cite{BDP2019} derived the CLT for the LSS of the large dimensional noncentral Fisher matrix. In this paper, we concentrate on
the outlier eigenvalues of the noncentral Fisher matrix defined in (\ref{noncentralF}). Specially, we will work with the following assumption,

	{\bf Assumption a}: Assume that $\boldsymbol{\Xi}$ is a $p\times n$ nonrandom matrix and the \textit{empirical spectral distribution} (ESD) of  $\boldsymbol{\Xi}\boldsymbol{\Xi}^{\ast}/n$ satisfies $H_n\overset{w}{\to}H$, ($w$ denoting weakly convergence), where $H$ is a non-random probability measure. In addition, the eigenvalues of $\boldsymbol{\Xi}\boldsymbol{\Xi}^{\ast}/n$ are subject to the condition
	\begin{align}\label{spikezhuan}
		l_{j_k+1}^{\boldsymbol{\Xi}}=l_{j_k+2}^{\boldsymbol{\Xi}}=\cdots=l_{j_k+m_k}^{\boldsymbol{\Xi}}={a_k},\quad k\in\{1,\cdots,K\}
	\end{align}
and ${a_k}$ satisfies the separation condition, that is,
	\begin{align}\label{sc}
		\mathop{min}\limits_{k\neq j}\left|\frac{a_k}{a_j}-1\right|>d,
	\end{align}
	where $d$ is positive constant and independent of $n$ and ${a_k}, k\in\{1,\cdots,K\}$ is allowed to grow at an order $o(\sqrt{n})$. In addition, $\mathcal{J}_k=\{j_k+1,\cdots,j_k+m_k\}$ denotes the set of ranks of $a_k$, where $m_k$ is the multiplicities of $a_k$ satisfying   $m_1+\cdots+m_K=M$, a fixed integer. 
\begin{remark}
	Note that ${a_k}, k\in\{1,\cdots,K\}$ can be located in any gap between the supports of $H$, which means that ${a_k}$ is not just the extreme eigenvalues of ${\boldsymbol{\Xi}\boldsymbol{\Xi}^{\ast}/n}$.
\end{remark}

The eigenvalues  $a_k, k\in\{1,\cdots,K\}$ are called as population spiked eigenvalues of the noncentral sample covariance matrix (\ref{noncentralS}) and the noncentral Fisher matrix (\ref{noncentralF}). We call these two matrices satisfying (\ref{spikezhuan})  the spiked noncentral sample covariance matrix and the spiked noncentral Fisher matrix, respectively. In fact, the spiked eigenvalues $a_k, k\in\{1,\cdots,K\}$ should have allowed to diverge at any rate, but under the restriction of our studying method, we have to assume them at a rate of $o(\sqrt{n})$. In this paper, we are devoted to exploring the properties of the limits of the sample spiked eigenvalues (corresponding to $a_k, k\in\{1,\cdots,K\}$) of the noncentral Fisher matrix. From now on, we call the sample spiked eigenvalues simply spiked eigenvalues when  no confusion can arise. 

To study the {\it principal component analysis} (PCA), \cite{Johnstone2001} proposed the spiked model based on the covariance matrix. The spiked model has been studied much further and extended in various random matrices such as Fisher matrix, sample canonical correlation matrix, separable covariance matrix, see \cite{Bao2018,DY2019} for more details. The main emphasis of the research of the spiked model is on the limits and the fluctuations of the spiked eigenvalues of these kinds of random matrices. \cite{Baik2005,Paul2007,BY2008,BY2012,BD2012,CaiXiaoPan(2017),JiangBai(2018)} focused on the spiked sample covariance matrix and \cite{WangYao2017,JHH2019,JHB2019} concentrated on the central spiked  Fisher matrices. 

The main contribution of the paper is the establishment of the limits and fluctuations of the spiked eigenvalues of the noncentral Fisher matrix under the Gaussian population assumption. Even better, we apply the above theoretical results to the {\it canonical correlation analysis} (CCA) and derive the limits and fluctuations of sample  
canonical correlation coefficients and give three consistent estimators, including the population spiked eigenvalues and the population canonical correlation coefficients. In addition, we study the properties of sample spiked eigenvalues of the noncentral sample covariance matrix, which should have its own interest.

The rest of the paper is organized as follows. In Section $2$, we define some notations and we present the LSD of some random matrices.  In Section $3$, we study the limits and fluctuations of spiked eigenvalues of the noncentral sample covariance matrix and noncentral Fisher matrix.  In Section $4$, we present the limits and fluctuations of spiked eigenvalues of sample canonical correlation matrix are investigated, give three estimators of the population spiked eigenvalues, {and conduct the actual data analysis about climate and geography by CCA}. {To show the correctness and rationality of theorems intuitively, we design a series of simulations in Section $5$.} In Section $6$, we summarize the main conclusions and the outlook. Section {$7$} presents technical proof. 

\section{Preliminaries}  \label{Pre}
In this section, we collect some notations and preliminary results or assumptions, which will be throughout the paper. Although some notations have been mentioned above, we still provide the precise  definitions here. 

\subsection{Basic notions}

For any $n\times n$ matrix $\mathbf{A}_n$ with only real eigenvalues, let $F_n$ be the \textit{empirical spectral distribution} (ESD) function of $\mathbf{A}_n$, that is,
$$F_n(x)=\frac{1}{n}Card\{i;\lambda_i^{\mathbf{A}_n}\leq x\},$$
where  $\lambda_i^{\mathbf{A}_n}$ denotes the $i$-th largest eigenvalue of $\mathbf{A}_n$.
If $F^{\mathbf{A}_n}$ has a limiting distribution $F$, then we call it the \textit{limiting spectral distribution} (LSD) of sequence $\{\mathbf{A}_n\}$.
For any function of bounded variation $G$ on the real line, its \textit{Stieltjes transform} (ST) is defined by
$$m(z)=\int\frac{1}{\lambda-z}dG(\lambda),~~z\in\mathbb{C}^{+}.$$
\subsection{Symbols and Assumptions}

In this paper, we study the spiked eigenvalues of the noncentral spiked sample covariance matrix $\mathbf C_n$ defined in (\ref{noncentralS}) and the noncentral spiked Fisher matrix $\mathbf F_p$ defined in (\ref{noncentralF}) with the matrix ${\mathbf \Xi}{\mathbf \Xi}^*/n$ satisfying (\ref{spikezhuan}). In order to distinguish the symbols of these three matrices clearly, we show the notations of eigenvalues and ST of the matrices in Table \ref{symbolmatrix}.

\begin{table*}[hp]
	\caption{Notations for $\boldsymbol{\Xi}\boldsymbol{\Xi}^*/n$, $\bbC_n$ and $\bbF_p$. }\label{symbolmatrix}
	\begin{tabular}{cccc}
		\hline
		Matrix & $\boldsymbol{\Xi}\boldsymbol{\Xi}^*/n$ &$\bbC_n$ &$\bbF_p$\\
		\hline
		LSD& $H$ & $F^{\bbC}$ &$F$\\
		\hline
		ST& $m_1$ & $m_2$ & $m_3$\\
		\hline
		Population spiked eigenvalue & $a_k$ & $a_k$  & $a_k$ \\
		\hline
		Sample eigenvalue & - & $l_i^{\bbC_n}$ & $l_i$\\
		\hline
		Limit & - & $\lambda_k^{\bbC}=\psi_{\bbC}(a_k)$ & $\lambda_k=\psi(a_k)=\psi_{\bbF}(\psi_{\bbC}(a_k))$\\
		\hline
	\end{tabular}
	
\end{table*}  
Throughout the paper,  we consider the following assumptions about the high-dimensional setting and the moment conditions. \\
{\bf Assumption b}: Assume that $p<n$, $p<N$ with $p/n=c_{1n}\rightarrow c_1\in(0,1), p/N=c_{2N}\rightarrow c_2\in(0,1)$, as $\min( p, n, N)\rightarrow \infty$.\\
{\bf Assumption c}: Assume that the matrix $\bbX_n$ and $\bbY_N$ are two independent arrays of independent standard Gaussian distribution variables $\{X_{jk}: 1\leq j\leq p, 1\leq k\leq n\}$ and $\{Y_{jk}: 1\leq j\leq p, 1\leq k\leq N\}$, respectively.  If $X_{ij}$ and $Y_{ij}$ are complex,
$E X_{ij}^2=0$, $E Y_{ij}^2=0$, $E |X_{ij}|^2=1$ and $E |Y_{ij}|^2=1$ are required.


\subsection{LSD for $\bbC_n$ and $\bbF_p$}\label{3.2}
To introduce the necessary conclusions and symbols in the following sections, we provide the LSD for $\bbC_n$ and $\bbF_p$ by the results in \cite{DS2007a,ZBY2015}. Note that similar results are obtained in \cite{BDP2019}.
According to Theorem $1.1$ of \cite{DS2007a}, the ST $m_2(z)$ of the LSD of $\bbC_n$ is the unique solution to the equation   
\begin{eqnarray}\label{m_2}
	m_2&=&\int\frac{dH(t)}{\frac{t}{1+c_1m_2}-(1+c_1m_2)z+1-c_1}\nonumber\\
	&=&\left(1+c_1m_2\right)m_1[(1+c_1m_2)((1+c_1m_2)z-(1-c_1))]
\end{eqnarray}
where $m_1(z)$ is the ST of $H$. For simplicity, we write $m_2(z)$ as $m_2$ in (\ref{m_2}). By Theorem $2.1$ of \cite{ZBY2015}, the ST $m_3(z)$ of the LSD of $\bbF_p$ satisfies the following equation:
\begin{eqnarray}\label{lsd3}
	m_3=\int\frac{dH(t)}{\frac{t}{1+(c_1+c_2z)m_3}+\frac{1-c_1}{1+c_2zm_3}-\frac{z(1+(c_1+c_2z)m_3)}{1+c_2zm_3}}.
\end{eqnarray}
From (6.12) in \cite{BS2010}, we have
\begin{align}\label{m_3}
	\frac{m_3(z)}{1+c_2zm_3(z)}=m_2\left(z(1+c_2zm_3(z))\right).
\end{align}

\section{Main Results}\label{mainresults}
In this section, we state our main results and briefly summarize our proof strategy. Our main results include the limits and the CLT of the spiked eigenvalues of the noncentral sample covariance matrix $\bbC_n$ and the noncentral Fisher matrix $\bbF_p$.
\subsection{Limits and fluctuations for the matrix $\bbC_n$}
The noncentral sample covariance matrix  $\bbC_n$ and its eigenvalues are arranged as a descending order as
\begin{align}\label{eigtuS1}
	l_1^{\bbC_n}\geq l_2^{\bbC_n} \geq \cdots \geq l_p^{\bbC_n}.
\end{align}

\begin{theorem}\label{limittuS1}
	If Assumption $[{\mathbf a}]-[{\mathbf c}]$ hold and the population spiked eigenvalues ${a_k}, \,k\in\{1,\cdots,K\}$ satisfies $\psi_{\bbC}'\left({a_k}\right)>0$, for $1\leq k\leq K$. Then we have 
	\begin{align}
		\frac{l_j^{\bbC_n}}{\psi_{\bbC}(a_k)}-1	\overset{a.s.}{\longrightarrow}0, \quad j\in\mathcal{J}_k\ 
	\end{align}	
	where
	\begin{align}\label{psituS1}
		{\psi_{\bbC}\left(a_k\right)}=a_k\left(1-c_{1n}\int\frac{1}{t-a_k}dH_n(t)\right)^2+(1-c_{1n})\left(1-c_{1n}\int\frac{1}{t-a_k}dH_n(t)\right),
	\end{align}
\end{theorem}
\begin{remark}
Considering that the convergence  $H_n\rightarrow H$ may be slow, in ${\psi_{\bbC}\left(a_k\right)}$, we use $H_n$ and $c_{1n}$ instead of $H$ and $c_1$, respectively. In following theorems related to the limits of spiked eigenvalues, we will take the same treatment. 
	\end{remark}

Having disposed of the limits of the sample spiked eigenvalues, we are now in a position to show the CLT for the sample spiked eigenvalues.
\begin{theorem}\label{CLTtuS1}
	Suppose the Assumption $[{\mathbf a}]-[{\mathbf c}]$ hold, the $m_k$-dimensional random vector
	\begin{equation*}
		\gamma_k^{\bbC_n}=\sqrt{n}\left\{\left(	\frac{l_j^{\bbC_n}}{\psi_{\mathbf C}\left(a_k\right)}-1\right)\frac{1}{\sqrt{\beta\theta_1}},\quad j\in\mathcal{J}_k\right\}
	\end{equation*}
converges weakly to the joint distribution of the $m_k$ eigenvalues of Gaussian random matrix $\boldsymbol\Omega$, where $\boldsymbol\Omega$ is a $m_k$-dimensional standard GOE (GUE) matrix.
If the samples are either real, $\beta=2$; or complex, $\beta=1$, and
	\begin{align*}
		\theta_1&=\frac{1}{[\lambda_k^{\bbC}\underline{m}_2'+\frac{a_k(1+c_1m_2+c_1\lambda_k^{\bbC}m_2')}{\lambda_k^{\bbC}(1+c_1m_2)^2}]^2}\\
		&\times\left(\underline{m}_2'+\frac{a_k^2c_1m_2'}{(\lambda_k^{\bbC})^2(1+c_1m_2)^4}+\frac{2a_k(1+\underline{m}_2+\lambda_k^{\bbC}\underline{m}_2')}{(\lambda_k^{\bbC})^2(1+c_1m_2)^2}\right),
	\end{align*}
	where $m_2'$ and $\underline{m}_2'$ are the  derivatives of $m_2$ and $\underline{m}_2$ at the point $\lambda_k^{\bbC}$, respectively, and 
	\begin{align}\label{munderline}
		\underline{m}_2(\lambda_k^{\bbC})=-\frac{1-c_1}{\lambda_k^{\bbC}}+c_1m_2(\lambda_k^{\bbC}).
	\end{align}
\end{theorem}
\begin{remark}
	{It is worth pointing out that  $\theta_1$ is equal to the (half of) variance of $\sqrt{n}(l_j^{\bbC_n}/\psi_{\bbC}(a_k)-1)$ when $a_k$ is a single eigenvalue. When $a_k$ is multiple, the limiting distribution of $\sqrt{n}(l_j^{\bbC_n}/\psi_{\bbC}(a_k)-1)$ is related to that of the eigenvalues of a GOE (GUE) matrix, the variance of whose diagonal elements is equal to $2\theta_1$ ($\theta_1$).} {In what follows, we also call $\theta_1$ as the scale parameter of the GOE (GUE) matrix.}
\end{remark}
\begin{remark}
\cite{D2020} and \cite{BDW2021}	focused on the noncentral spiked sample covariance matrix. \cite{D2020} studied the limits and rates for the spiked eigenvalues and vectors of the noncentral sample covariance matrix. \cite{BDW2021} was devoted to the fluctuation of the spiked vectors of
the noncentral sample covariance matrix. Note that both of the works are under the condition of finite rank, specially, the rank of the ${\mathbf \Xi}$ is finite.  Compared with \cite{D2020} and \cite{BDW2021},  our assumptions about ${\mathbf \Xi}$ is more general.
\end{remark}

\subsection{Limits and fluctuations for the noncentral Fisher matrix $\bbF_p$}

Having disposed of the noncentral  spiked sample covariance matrix ${\mathbf C}_n$,  we can now return to the noncentral Fisher matrix $\bbF_p$.  The eigenvalues of noncentral Fisher matrix $\bbF_p$ are sorted in descending order as
\begin{eqnarray}
	l_1\geq l_2\geq\cdots\geq l_p.
\end{eqnarray}

\begin{theorem}\label{limitsF}
	Let the Assumption $[{\mathbf a}]-[{\mathbf c}]$ hold, the noncentral Fisher matrix $\bbF_p$ is defined in (\ref{noncentralF}). If $a_k$ satisfies $\psi_{\mathbf F}'(\psi_{\mathbf C}(a_k))>0$ and $\psi_{\mathbf C}'(a_k)>0$, for $1\leq k\leq K$, then we have 
	\begin{align*}
		\frac{l_{j}}{\psi_{\mathbf F}\left(\psi_{\mathbf C}(a_k)\right)}-1\overset{a.s.}{\longrightarrow} 0, \quad j\in\mathcal{J}_k 
	\end{align*}
	where 
	\begin{eqnarray}\label{psiF}
		\psi_{\mathbf F}\left(x\right)=\frac{x}{1+c_{2n}\cdot x\cdot m_{2}^0\left(x\right)},
	\end{eqnarray}
	$\psi_{\mathbf F}'(\cdot)$ is the derivative of $\psi_{\boldsymbol{F}}(\cdot)$, $m_2^0$ and $m_{2}$ are same, but $c_1$, $c_2$ and $H$ replaced by $c_{1n}$, $c_{2n}$ and $H_n$, respectively.
\end{theorem}
The task is now to show the CLT for the sample spiked eigenvalues of the noncentral Fisher matrix $\bbF_p$.
\begin{theorem}\label{CLTnonF}
	Suppose that the Assumption $[{\mathbf a}]-[{\mathbf c}]$ hold, the $m_k$ dimensional random vector
	\begin{equation*}
		\gamma_{k}^{\bbF_p}\overset{\triangle}{=}\sqrt{{n}}\left\{\left(\frac{l_j-\psi_{\mathbf F}(\psi_{\boldsymbol C}(a_k))}{\psi_{\mathbf F}(\psi_{\boldsymbol C}(a_k))}\right)\frac{1}{\sqrt{\beta\theta_2}}, \quad j\in\mathcal{J}_k\right\},
	\end{equation*}
	converges weakly to the joint distribution of the eigenvalues of Gaussian random matrix $\boldsymbol{\Omega}$ where 
	\begin{equation}\label{theta}
		\theta_2=\frac{c_2}{c_1\cdot\vartheta}+\left[\frac{1-c_2(\lambda_{k}^{\bbC})^2m_2'(\lambda_{k}^{\bbC})}{1+c_2\lambda_{k}m_3(\lambda_{k})}\right]^2\theta_1,
	\end{equation}
$\theta_1$ is defined in Theorem \ref{CLTtuS1}
 and $\vartheta$ satisfies
\begin{align}
\vartheta&=1+2\lambda_kc_2m_3(\lambda_k)+c_2\lambda_k^2m_3'(\lambda_k)\label{theta_2}
\end{align}
\end{theorem}

\section{Applications}
In this section, we discuss some  applications of our results in limiting properties of sample canonical correlations coeffcients, estimators of the population spiked eigenvalues and the population canonical correlations coeffcients. At the end, we present an experiment on a real environmental variables for world countries data.
\subsection{Limits and fluctuations for the sample canonical correlation matrix}

The CCA is the general and favorable method to investigate the relationship between two random vectors. Under the high-dimensional setting and the Gaussian assumption, \cite{Bao2018} studied the limiting properties of sample canonical correlation coefficients. Under the sharp moment condition,   \cite{Yang2020,Yang2021} prove the largest eigenvalue of sample canonical correlation matrix converges to T-W law to test the independence of  random vectors with two different structures. {Moreover, \cite{Yang21Limiting} shows the limiting distribution of the spiked eigenvalues depends on the fourth cumulants of the population distribution.}
 Note that \cite{BHHJZ2020,Yang2020,Yang2021,Yang21Limiting} only focused on the finite rank case, where the number of the positive population canonical correlation coefficients of two groups of high-dimensional Gaussian vectors are finite. However, we popularize finite rank case to infinite rank case, in other words, we get the limits and fluctuations of the sample canonical correlation coefficients under the infinite rank case. In the following, we will introduce the application about limiting properties of sample canonical correlations coeffcients in great detail. 

Let ${\mathbf{z}_i}=(\mathbf x_i^T, \mathbf y_i^T)^T, i=1,\cdots, n$, be independent observations from  a $(p+q)$-dimensional Gaussian distribution 
with mean zero and covariance matrix
\begin{align*}\label{popcov}
	\boldsymbol \Sigma=\left(
	\begin{array}{cc}
		\boldsymbol \Sigma_{xx} & \boldsymbol \Sigma_{xy} \\
		\boldsymbol \Sigma_{yx} & \boldsymbol \Sigma_{yy} \\
	\end{array}
	\right),
\end{align*}
where ${\mathbf x}_i$ and ${\mathbf y}_i$ are $p$-dimensional and $q$-dimensional vectors with the population covariance matrices $\boldsymbol \Sigma_{xx}$ and $\boldsymbol \Sigma_{yy}$, respectively. Without loss of generality, we assume that  $p\le q$.
Define the corresponding sample covariance matrix as

\begin{align}
{\boldsymbol S}_n=\frac{1}{n}\sum_{i=1}^{n}{\boldsymbol z}_i{\boldsymbol z}_i^T,
\end{align}
which can be formed as
\begin{align*}
	{\boldsymbol S}_n=\left(
	\begin{array}{cc}
		{\boldsymbol S}_{xx} & {\boldsymbol S}_{xy} \\
		{\boldsymbol S}_{yx} & {\boldsymbol S}_{yy} \\
	\end{array}
	\right)
	=\frac{1}{n}\left(
	\begin{array}{cc}
		{\boldsymbol X}{\boldsymbol X}^T & {\boldsymbol X}{\boldsymbol Y}^T \\
		{\boldsymbol Y}{\boldsymbol X}^T & {\boldsymbol Y}{\boldsymbol Y}^T \\
	\end{array}
	\right)
\end{align*}
with
\begin{align*}
	{\boldsymbol X}=({\mathbf x}_1,\cdots, {\mathbf x}_n)_{p\times n}, \quad  {\boldsymbol Y}=({\mathbf y}_1,\cdots, {\mathbf y}_n)_{q\times n},
\end{align*}In the sequel, 
${\boldsymbol \Sigma}_{xx}^{-1}{\boldsymbol \Sigma}_{xy}{\boldsymbol \Sigma}_{yy}^{-1}{\boldsymbol \Sigma}_{yx}$ is called as the population canonical correlation matrix and its eigenvalues  are denoted by \begin{align}\label{pcca}
	1>\rho_1^2\geq \rho_2^2 \geq \cdots \geq \rho_p^2.
\end{align}
By the singular value decomposition, we have that
\begin{align}\label{P1P2}
\boldsymbol \Sigma_{xx}^{-\frac{1}{2}}\boldsymbol \Sigma_{xy}\boldsymbol \Sigma_{yy}^{-\frac{1}{2}}={\mathbf P}_1{\boldsymbol \Lambda}{\mathbf P}_2^T
\end{align}
where
\begin{align}\label{Lam}
{\boldsymbol \Lambda }=\left(
\begin{array}{cc}
{\boldsymbol \Lambda }_{11}&  {\mathbf 0}_{12}\\
\end{array}
\right),
\end{align}
${\boldsymbol \Lambda }_{11}= \text{diag}(\rho_1,\rho_2, \cdots, \rho_p)$, ${\mathbf 0}_{12}$ is a $p\times (q-p)$ zero matrix, ${\mathbf {P}}_1$ and ${\mathbf {P}}_2$ are orthogonal matrix with size $p\times p$ and $q \times q$, respectively. It follows that $\rho^2_1, \rho^2_2, \cdots, \rho^2_p$ are also the eigenvalues of the diagonal matrix ${\boldsymbol \Lambda}{\boldsymbol \Lambda}^T$. According to Theorem 12.2.1 of \cite{Anderson2003}, the nonnegative square roots $\rho_1,\cdots, \rho_p$ are the population canonical correlation coefficients. Correspondingly, 
${\mathbf S}_{xx}^{-1}{\mathbf S}_{xy}{\mathbf S}_{yy}^{-1}{\mathbf S}_{yx}$  is called as the sample canonical correlation matrix and its eigenvalues are denoted by
\begin{align}\label{scca}
	\lambda_1^2\geq \lambda_2^2 \geq \cdots \geq \lambda_p^2.\end{align}
The following theorem describes the function relation between sample canonical correlation coefficients and the eigenvalues of a special noncentral Fisher matrix. 

\begin{theorem}\label{translation}[Theorem 1 in \cite{BHHJZ2020}]
	Suppose that $\lambda_i^2,  i=1,\cdots, p$, is the ordered eigenvalue of the sample canonical correlation matrix ${\mathbf S}_{xx}^{-1}{\mathbf S}_{xy}{\mathbf S}_{yy}^{-1}{\mathbf S}_{yx}$. Then, there exists a noncentral Fisher matrix $\boldsymbol F({\mathbf \Xi})$ whose eigenvalue $l_i$ satisfies
	$l_i=g(\lambda_i)\overset{\triangle}{=}\frac{(n-q)\lambda_i^2}{q(1-\lambda_i^2)}, ~i=1,\cdots, p$
	and ${\mathbf \Xi}$ is the noncentral parameter matrix,
	\begin{align}\label{Noncenp}
		{\mathbf \Xi}=\frac{n}{q}\mathbf{T}(n^{-1}\widehat{\mathbf Y}\widehat{\mathbf Y}^T)\mathbf{T}^T
	\end{align}
	where ${\mathbf T}{\mathbf T}^{T} =diag\left(\rho^2_1/(1-\rho^2_1),\rho^2_2/(1-\rho^2_2),\cdots,\rho^2_p/(1-\rho^2_p)\right)$ and $\widehat{\mathbf Y}$ is a $p\times n$ matrix and contains i.i.d. elements with standard Gaussian distribution.  
\end{theorem}

Combining Theorem \ref{translation} with the properties of the noncentral Fisher matrix, we obtain the limits and fluctuations of the sample canonical correlation coefficients under the following assumption.

{\bf Assumption d}: {The \textit{empirical spectral distribution} (ESD) of ${\boldsymbol \Sigma}_{xx}^{-1}{\boldsymbol \Sigma}_{xy}{\boldsymbol \Sigma}_{yy}^{-1}{\boldsymbol \Sigma}_{yx}$,  $\mathcal H_n$, tends to proper probability measure $\mathcal H$,  if  $\min( p, q, n)\rightarrow \infty$. Assume that  $\rho^2_i,~i=1,\dots,p$ is subject to the condition
	\begin{align}\label{ak}
		\alpha_k=\rho_{m_{k-1}+1}^2=\cdots=\rho_{m_{k-1}+m_k}^2, \quad k\in\{1,\cdots,K\},
	\end{align}
	where $\alpha_k$ is out of the support of $\mathcal H$ and satisfies the separation condition defined in (\ref{sc}). $M=\sum_{i=1}^{K}m_i$ is a fixed positive integer with convention $m_0=0$. 
	In addition, $\alpha_1$ is allowed to be with the order $1-o(n^{-1/2})$.}

\begin{remark}
	{Note that the noncentral parameter matrix (\ref{Noncenp}) is random, so the assumption about multiple roots in (\ref{ak}) is not reasonable. In the following Assumption d$'$,  we replace the comndition of multiple roots by single root. However, we think the results of the limits and fluctuations for the multiple roots case are correct and our guess will be verifed by some simulations in the following section.}
\end{remark}

{\bf Assumption d$'$}:The assumptions are same as Assumption d, the additional assumption  $m_i=1, \quad i\in\{1,\cdots,K\}$.
\begin{theorem}\label{limitsCCA}
	Under the conditions stated in Theorem \ref{translation}, if moreover Assumption d' holds, and $\alpha_k$ satisfies $\psi_{\boldsymbol{\Xi}}'(f(\alpha_k))>0$, $\Psi_{\mathbf C}'(\alpha_k)>0$, and $\Psi'(\alpha_k)>0$ for $1\leq k\leq K$, then we have 
	\begin{align}
		\frac{\lambda_i^2}{t(\alpha_k)}-1\overset{a.s.}{\longrightarrow} 0, \quad i\in{\mathcal J}_k 
	\end{align}	
	where	 
	\begin{align*}\label{Psi}
		&t(x)=g^{-1} \circ \Psi(x),\quad \Psi(x)=\psi_{\boldsymbol F}\circ\psi_{{\boldsymbol C}}\circ\psi_{\boldsymbol \Xi}\circ f(x),\\
		&\Psi_{\bbC}(x)=\psi_{\bbC}\circ\psi_{\boldsymbol{\Xi}}\circ f(x), \quad f(x)=\frac{n}{q}\frac{x}{1-x}, \quad\psi_{\boldsymbol F}(x)=\frac{x}{1+\frac{p}{n-q}\cdot x\cdot m_{\boldsymbol C}(x)};\nonumber \\
		&\psi_{\boldsymbol \Xi}(x)=x\left(1+\frac{p}{n}\int\frac{t}{x-t}d\widetilde{\mathcal H}(t)\right), \quad \widetilde{\mathcal H}(x)=\mathcal H_n\left(\frac{qx}{n+qx}\right),\nonumber\\
		&\psi_{\boldsymbol C}(x)=x\left(1\!-\!\frac{p}{q}\int\frac{1}{t-x}dF^{p/n, {\tilde{\mathcal{H}}}}_{mp}(t)\right)^2\!-\!\left(1\!-\!\frac{p}{q}\right)\left(1\!-\!\frac{p}{q}\int\frac{1}{t-x}dF^{p/n, {\tilde{\mathcal{H}}}}_{mp}(t)\right),\nonumber
	\end{align*}
	and $g^{-1} \circ \Psi$ stands for the composition of $\Psi(\cdot)$ and  $g^{-1}(\cdot)$. $F_{mp}^{p/n, {\tilde{\mathcal{H}}}}$ denotes  M-P law with  the parameter $p/n$ and ${\tilde{\mathcal{H}}}$. Moreover, $m_{\boldsymbol C}(\cdot)$ stands for the unique solution (\ref{m_2}) with $c_1$, $H$ replaced by $p/q$, $F_{mp}^{p/n, {\tilde{\mathcal{H}}}}$, respectively. 
\end{theorem}

It is worth noting that we can conclude the results in Theorem 1.8 of \cite{Bao2018} by Theorem \ref{limitsCCA} or Theorem 3.1 in \cite{WangYao2017}, when the LSD $\mathcal{H}$ of ${\boldsymbol \Sigma}_{xx}^{-1}{\boldsymbol \Sigma}_{xy}{\boldsymbol \Sigma}_{yy}^{-1}{\boldsymbol \Sigma}_{yx}$ degenerates to $\delta_{\{0\}}$.
\begin{corollary}\label{limitsCCAcor}
	Let the Assumption [${\mathbf c}$] and [${\mathbf d'}$] hold, furthermore, the LSD $\mathcal H$ degenerates to $\delta_{\{0\}}$ and $\alpha_k (k=1, \cdots, K)$ satisfies $ \alpha_k>\alpha_r,$
	then we have
	\begin{align}
	\frac{\lambda_i^2}{\phi(\alpha_k)}-1\overset{a.s.}{\longrightarrow}0,\quad i\in{\mathcal J}_k,
	\end{align}
where  
	\begin{align*}
		\phi(\alpha_k)&=\frac{[\alpha_k(1-r_1)+r_1][\alpha_k(1-r_2)+r_2]}{\alpha_k}\\ \alpha_r&=\sqrt{\frac{r_1r_2}{(1-r_1)(1-r_2)}}\\
		r_1&=p/n,\quad r_2=q/n.
	\end{align*}
\end{corollary}
Having disposed of the results of limits of the square of sample canonical correlation coefficients $\lambda_i^2,  i\in{\mathcal J}_k$ associated to the $\alpha_k (k=1, \cdots, K)$, we can now return to show the CLT for them.
\begin{theorem}\label{cltcca}
	Let the Assumption $[{\mathbf c}]$ and $[{\mathbf d'}]$ hold, and $\lambda_i^2$ is the square of eigenvalues of sample canonical correlation matrix. We set {$p/q\to c_3$} and {$p/(n-q)\to c_4$} as $n$ tends to infinity, then the $m_k$-dimensional random vector
	\begin{align*}
		\gamma_k=\sqrt{q}\left\{\left(\frac{\lambda_i^2-t(\alpha_k)}{t(\alpha_k)}\right)\frac{1}{\sqrt{\beta\eta}},\quad i\in\mathcal{J}_k\right\}
	\end{align*}
	converges weakly to the joint distribution of the eigenvalues of Gaussian random matrix $\boldsymbol{\Omega}$, and
		\begin{align*}
			&\eta=\left[\frac{c_4}{(c_3\cdot\eta_2)}+\left(\frac{1-c_4(\Psi_{\bbC}(\alpha_k))^2m_{\bbC}'(\Psi_{\bbC}(\alpha_k))}{1+c_4\Psi(\alpha_k)m_{\bbF}(\Psi(\alpha_k))}\right)^2\eta_1+\eta_3\right]\\
			&\quad\quad\times \frac{c_3^2c_4^2{\Psi(\alpha_k)^2}}{\left[c_3+c_4\Psi(\alpha_k)\right]^4{[t(\alpha_k)]^2}} ,\\ 
			&\eta_1=\left[\Psi_{\bbC}(\alpha_k)\underline{m}_{\bbC}'+\frac{\psi_{\boldsymbol{\Xi}}(f(\alpha_k))(1+c_3m_{\bbC}+c_3\Psi_{\bbC}(\alpha_k)m_{\bbC}')}{\Psi_{\bbC}(\alpha_k)(1+c_3m_{\bbC})^2}\right]^{-2} \nonumber \\
			&\times\left(\underline{m}_{\bbC}'\!+\!\frac{(\psi_{\boldsymbol{\Xi}}(f(\alpha_k)))^2c_3m_{\bbC}'}{(\Psi_{\bbC}(\alpha_k))^2(1\!+\!c_3m_{\bbC})^4}\!+\!\frac{2\psi_{\boldsymbol{\Xi}}(f(\alpha_k))(1\!+\!\underline{m}_{\bbC}\!+\!\Psi_{\bbC}(\alpha_k)\underline{m}_{\bbC}')}{(\Psi_{\bbC}(\alpha_k))^2(1\!+\!c_3m_{\bbC})^2}\right),\\
			&\eta_2=1+2\Psi(\alpha_k)c_4m_{\bbF}(\Psi(\alpha_k))+c_4(\Psi(\alpha_k))^2m_{\bbF}'(\Psi(\alpha_k)),\\
			&\eta_3=\frac{\frac{q}{n}(\psi_{\bbF}'(\Psi_{\bbC}(\alpha_k)))^2(\psi_{\bbC}'(\psi_{\boldsymbol{\Xi}}(f(\alpha_k))))^2}{\psi_{\boldsymbol{\Xi}}^2(f(\alpha_k))\underline{m}'(\psi_{\boldsymbol{\Xi}}(f(\alpha_k)))}\frac{\psi_{\boldsymbol{\Xi}}^2(f(\alpha_k))}{\Psi^2(\alpha_k)}
		\end{align*}
		where $m_{\bbC}$ and $m_{\bbC}'$ denote the value and derivative at the point $\Psi_{\bbC}(\alpha_k)$, respectively, $m_{\bbF}$ stands for the unique solution (\ref{lsd3}) with $c_1$, $c_2$, $H$ replaced by $p/q$, $p/(n-q)$, $F_{mp}^{p/n, H}$, respectively.
\end{theorem}

{The proof of this theorem can be drawn by the Delta method, Theorem \ref{CLTnonF} and Theorem \ref{translation} and it will be postponed in subsection \ref{CLT4}.}

\subsection{Estimators of the population spiked eigenvalues}
In this section, we develop two consistent estimators of population spiked eigenvalues $a_k$ defined in (\ref{ak}), which are derived by the results of the noncentral sample covariance matrix and the noncentral Fisher matrix, respectively. 
At first, we proceed to show the estimator of population spiked eigenvalues for the noncentral sample covariance matrix. From the conclusion of Theorem \ref{limitsF}, we have 
\begin{align*}
	\lambda_i^{\bbC}=a_k\left(1-c_1m_1(a_k)\right)^2+\left(1-c_1\right)\left(1-c_1m_1(a_k)\right),
\end{align*}
and
\begin{align*}
	a_k=\lambda^{\bbC}_k\left(1+c_1m_2(\lambda^{\bbC}_k)\right)^2-(1-c_1)\left(1+c_1m_2(\lambda^{\bbC}_k)\right).
\end{align*}
It is sufficient to consider the estimators of $\lambda^{\bbC}_k$ and $m_2(\lambda^{\bbC}_k)$, which are denoted by $\hat{\lambda}^{\bbC}_k$ and $\hat{m}_2(\hat{\lambda}^{\bbC}_k)$, respectively. We adopt an approach similar to that in \cite{JiangBai(2018)} to estimate $\hat{m}_2(\hat{\lambda}^{\bbC}_k)$.
Define $r_{ik}=|\hat\lambda_{i}^{\bbC}-\hat\lambda_{k}^{\bbC}|/|\hat\lambda_{k}^{\bbC}|$ and the set $\mathcal {J}_k=\{i\in (1,\cdots, p): r_{ik}\leq 0.2 \}$ and $\tilde c_1=(p-|\mathcal {J}_k|)/n$; then,
	\begin{align}\label{xstieljes}
		\hat m_2(\hat\lambda_{k}^{\bbC})=\frac{1}{p-|\mathcal {J}_k|}\sum\limits_{i\notin \mathcal J_k}(\hat\lambda_{i}^{\bbC}-\hat\lambda_{k}^{\bbC})^{-1}
	\end{align}
	is a good estimator of $m_2(\lambda_k^{\bbC})$,
	where the set $\mathcal {J}_k$ is  selected to avoid  the effect of multiple roots and to make the estimator more accurate. Then we have
	\begin{align}\label{estimatorfornonS}
		\hat{a}_k=\hat{\lambda}^{\bbC}_k\left(1+\tilde{c}_1\hat{m}_2(\hat{\lambda}^{\bbC}_k)\right)^2-(1-\tilde{c}_1)\left(1+\tilde{c}_1\hat{m}_2(\hat{\lambda}^{\bbC}_k)\right). 
\end{align}

We now turn to give the other estimator for the population spiked eigenvalues by the results of the noncentral Fisher matrix. From the conclusion of Theorem \ref{limitsF}, we have 
\begin{align*}
	\lambda_i=\frac{\psi_{\bbC}(a_k)}{1+c_2\psi_{\bbC}(a_k)m_2(\psi_{\bbC}(a_k))},
\end{align*}
then,
\begin{align}\label{41}
	\psi_{\bbC}(a_k)=\lambda_i(1+c_2\lambda_im_3(\lambda_i)).
\end{align}
According to Theorem \ref{limittuS1}, we know 
\begin{align*}
	\psi_{\bbC}(a_k)=a_k(1-c_1m_1(a_k))^2+(1-c_1)(1-c_1m_1(a_k)),
\end{align*}
then,
\begin{align*}
	a_k=\psi_{\bbC}(a_k)\left[1+c_1m_2\psi_{\bbC}(a_k)\right]^2-(1-c_1)\left[1+c_1m_2\psi_{\bbC}(a_k)\right].
\end{align*}
Note that (\ref{41}) and $\hat{\lambda}_k$ is the natural estimator of $\lambda_k$, we set $\tilde{a}_k$ as the estimator of $\psi_{\bbC}(a_k)$ and have 
	\begin{align*}
		\tilde{a}_k=\hat{\lambda}_k(1+\tilde c_2\hat{\lambda}_k\hat{m}_3(\hat{\lambda}_k)),
	\end{align*}
	where $\hat{m}_3(\hat\lambda_k)$ are the estimators of $m_3(\lambda_k)$. Similar considerations in (\ref{xstieljes}) are applied here, we have $r_{ik}=|\hat\lambda_{i}-\hat\lambda_{k}|/|\hat\lambda_{k}|$ and the set $\mathcal {J}_k=\{i\in (1,\cdots, p): r_{ik}\leq 0.2\}$ and $\tilde c_2=(p-|\mathcal {J}_k|)/N$; then,
	\begin{align}\label{x3stieljes}
		\hat m_3(\hat\lambda_{k})=\frac{1}{p-|\mathcal {J}_k|}\sum\limits_{i\notin \mathcal J_k}(\hat\lambda_{i}-\hat\lambda_{k})^{-1}
	\end{align}
	is a good estimator of $m_3(\lambda_{k})$. 
	Then we have,
	\begin{align}\label{estimatorfornonF}
		\hat{a}_k=[\tilde{a}_k(1+c_1\hat{m}_2(\tilde{a}_k))+(1-c_1)](1+c_1\hat{m}_2(\tilde{a}_k)),
	\end{align}
	where $\hat{m}_2(\tilde{a}_k)$ is the estimator of $m_2$ at $\psi_{\bbC}(a_k)$, by (\ref{m_3}) we have
	\begin{align}\label{x2stieljes}
		\hat m_2(\tilde{a}_k)=\hat m_3/(1+\tilde{c}_2\hat\lambda_k\hat m_3).
	\end{align}


\subsection{{The} environmental variables for world countries data}
To illustrate the application of canonical correlation, we apply our result to the environmental variables for world countries data\footnote{The data can be download from: https://www.kaggle.com/zanderventer/environmental-variables-for-world-countries.}. Deleting the samples missing value and the variables related to others, then, we get $188$ samples with $18$ variables.  We divide the $18$ variables into two groups, one ($\mathbf x$, $p=7$) contains the elevation above sea level, percentage of country covered by cropland, percentage cover by trees $>$ $5$m in height and so on, the other one ($\mathbf y$, $q=11$) contains annual precipitation, temperature mean annual, mean wind speed, average cloudy days per year and so on. We want to explore the relationship between the geographical conditions $\mathbf x$ and climatic conditions $\mathbf y$, so we test their first canonical correlation coefficient being zero or not, i.e.,
\begin{align*}
H_0 : \rho_1^2=0\quad \text{v.s.}\quad H_1 : \rho_1^2>0.
\end{align*}
Therefore, we take the largest eigenvalue $\lambda_1^2$ of the sample canonical correlation matrix as the test statistic. By the formula (2.1) in \cite{Bao2018}, the normalized largest eigenvalues of the sample CCA matrix tends to the T-W law under the null hypothesis. The eigenvalues of ${\mathbf S}_{xx}^{-1}{\mathbf S}_{xy}{\mathbf S}_{yy}^{-1}{\mathbf S}_{yx}$ are as follows:

\begin{table}[h]
	\centering
	\caption{The sample canonical correlation coefficients of the real data}
	\begin{tabular}{p{1cm}p{1cm}p{1cm}p{1cm}p{1cm}p{1cm}p{1cm}}
		\toprule  
		$\lambda_1^2$ & $\lambda_2^2$ &$\lambda_3^2$& $\lambda_4^2$ & $\lambda_5^2$ & $\lambda_6^2$ & $\lambda_7^2$ \\
		\midrule  
		{0.9152} & {0.7755}& {0.4560} &{0.4034}& {0.2548} & {0.2247} & {0.0492} \\
		\bottomrule 
	\end{tabular}
\end{table}

According to the data, we obtain that $p$-value {approaches} zero. Thus we have strong evidence to reject the null hypothesis and conclude that the geographical conditions relate to climatic conditions. {Moreover, we use Algorithm \ref{algorithm1} to give an estimator  of the population canonical correlation coefficients,
\begin{align}\label{estimatorforCCA}
\hat{\rho}_i^2=\frac{q/n*\hat{a}_i}{1+q/n*\hat{a}_i},\quad i=1,\cdots,p.
\end{align}	
 By (\ref{estimatorforCCA}), we get the population canonical correlation coefficients $\hat{\rho}_1^2=0.9064$, which implies that the correlation between geographical conditions and climatic conditions is  strong.}

\begin{algorithm}[h]
	\caption{estimate population canonical correlation coefficients}\label{algorithm1} 
	\hspace*{0.02in} {\bf Input:} 
	the samples ${\mathbf x_i,\mathbf y_i, i=1,\cdots,n}$\\
	\hspace*{0.02in} {\bf Output:} 
	$\hat{\rho}_i^2, i=1,\cdots,p$
	\begin{algorithmic}[1]
		\State $\{l_p^{SCC}\leq\cdots\leq l_1^{SCC}\}$=Eigenvalue(${\mathbf S}_{xx}^{-1}{\mathbf S}_{xy}{\mathbf S}_{yy}^{-1}{\mathbf S}_{yx}$);
		\State $l_i^{F}=l_i^{SCC}/(1-l_i^{SCC})*(n-q)/q\quad i=1,\cdots,p;$
		\State Use step 2 and (\ref{estimatorfornonF}) to get the estimator for the population spiked eigenvalues of the noncentral Fisher matrix
		\State $\hat{\rho}_i^2=\frac{q/n*\hat{a}_i}{1+q/n*\hat{a}_i},\quad i=1,\cdots,p$
		\State \Return result $\hat{\rho}_i^2$
	\end{algorithmic}
\end{algorithm}

\section{Simulation}
We conduct simulations that support the theoretical results and illustrate the accuracy of the estimators. The simulations are divided into two sections, one is to verify the accuracy of Theorem \ref{CLTtuS1}, \ref{CLTnonF} and  \ref{cltcca}, and another is confirm the performance of the estimators in (\ref{estimatorfornonS}), (\ref{estimatorfornonF}) and (\ref{estimatorforCCA}).
\subsection{Simulations for the asymptotic normality}\label{Snormal}

 In this section, we assume the  eigenvalues of ${\mathbf \Xi}{\mathbf \Xi}^*/n$ satisfy 
\begin{align}\label{Sim1model}
10,\, 7.5,\, \underbrace{1,\,\cdots,\,1}_{p-2},
\end{align}	
in this situation, we set $\boldsymbol{\Lambda}$ defined in ({\ref{Lam}}) satisfying
\begin{align}\label{Sim11model}
\boldsymbol{\Lambda}=(\sqrt{10/11},\sqrt{15/17},\sqrt{1/2},\cdots,\sqrt{1/2}).
\end{align}

\begin{figure}[!h]	
	\centering
	\includegraphics[width=2.0in]{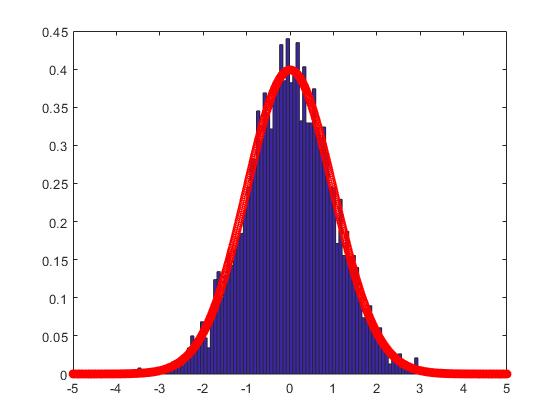}
	\includegraphics[width=2.0in]{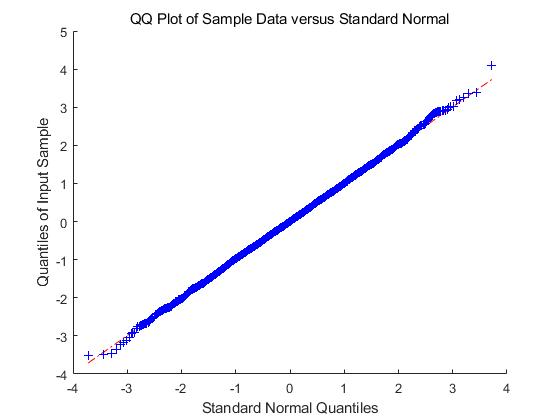}\\
	\caption{{The asymptotic normality of the largest eigenvalue of noncentral sample covariance matrix with $(p,n)=(200,2000)$.}}
	\label{penG1}
\end{figure}

\begin{figure}[!h]	
	\centering
	\includegraphics[width=2.0in]{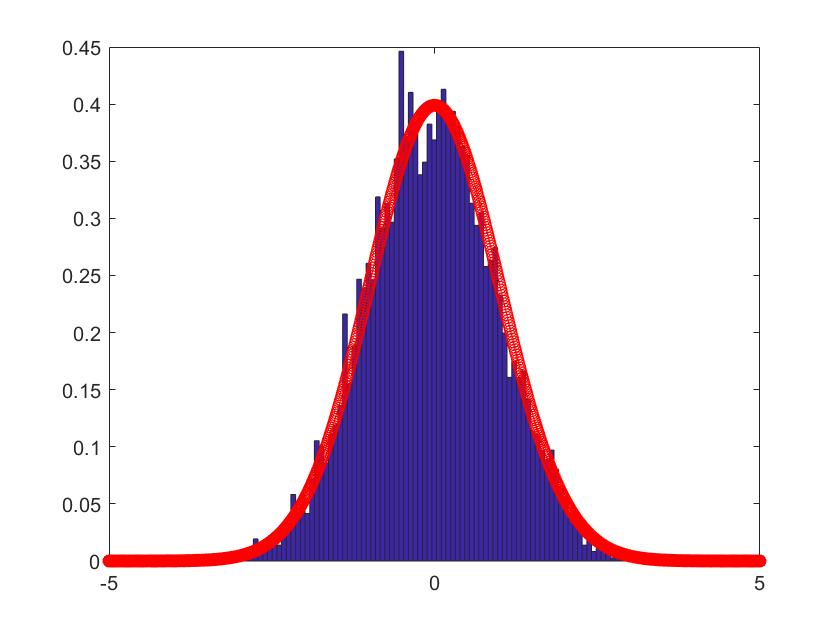}
	\includegraphics[width=2.0in]{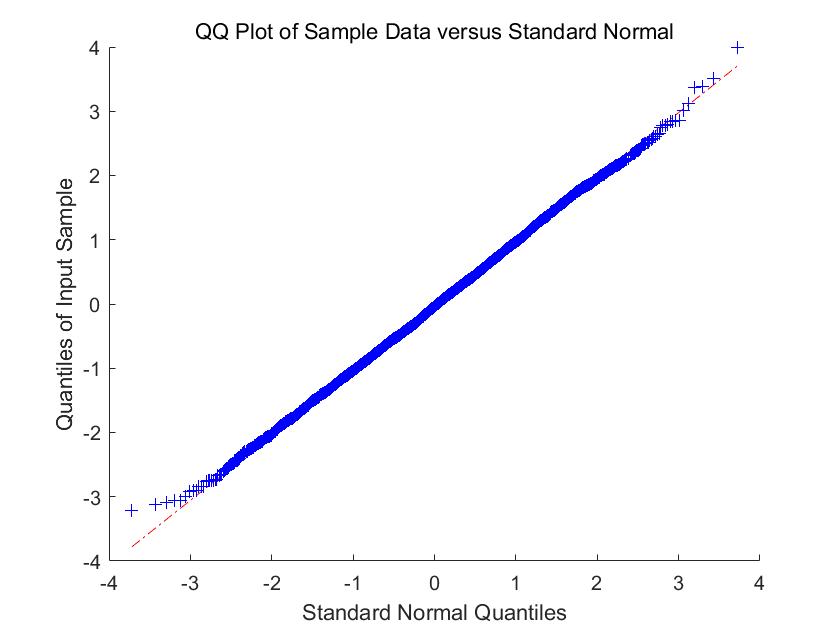}
	\caption{{The asymptotic normality of the second largest eigenvalue of the noncentral sample covariance matrix with $(p,n)=(200,2000)$.}}
	\label{penG2}
\end{figure}

\begin{figure}[!h]	
	\centering
	\includegraphics[width=2.0in]{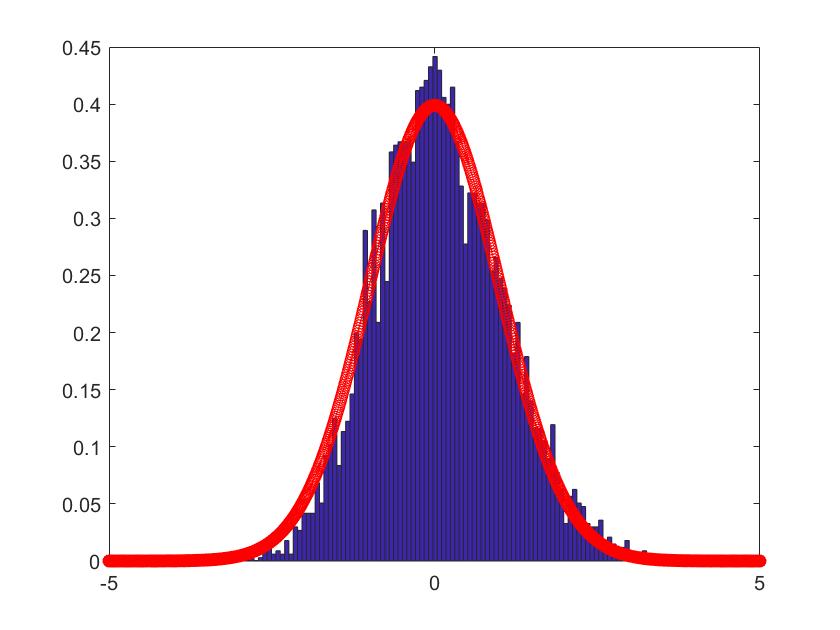}
	\includegraphics[width=2.0in]{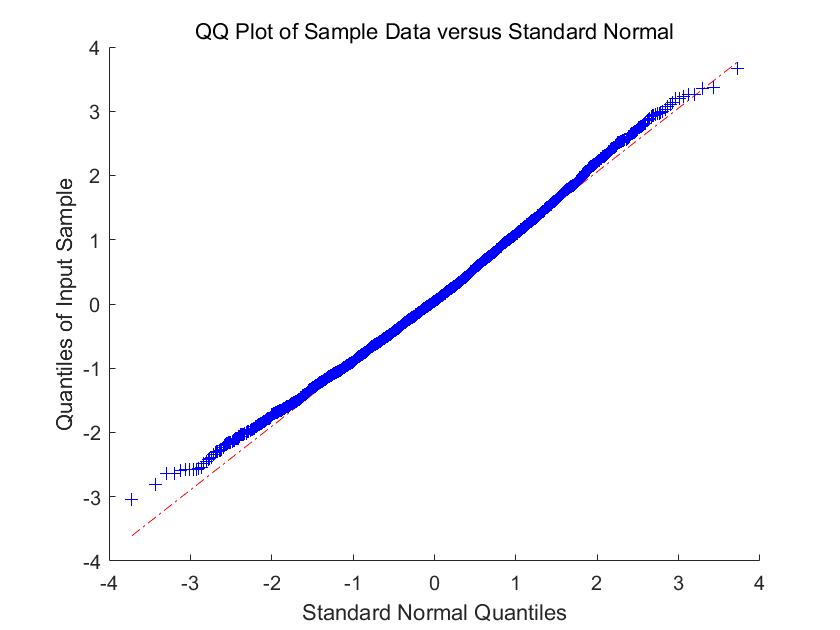}\\
	\caption{{The asymptotic normality of the largest eigenvalue of noncentral Fisher matrix with $(p,n,N)=(200,2000,1000)$.}}
	\label{penG3}
\end{figure}

\begin{figure}[!h]	
	\centering
	\includegraphics[width=2.0in]{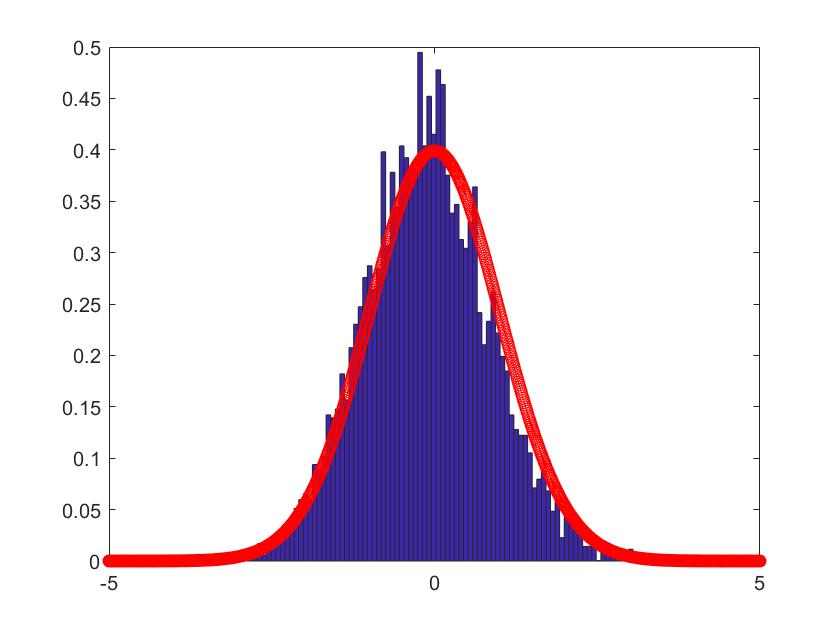}
	\includegraphics[width=2.0in]{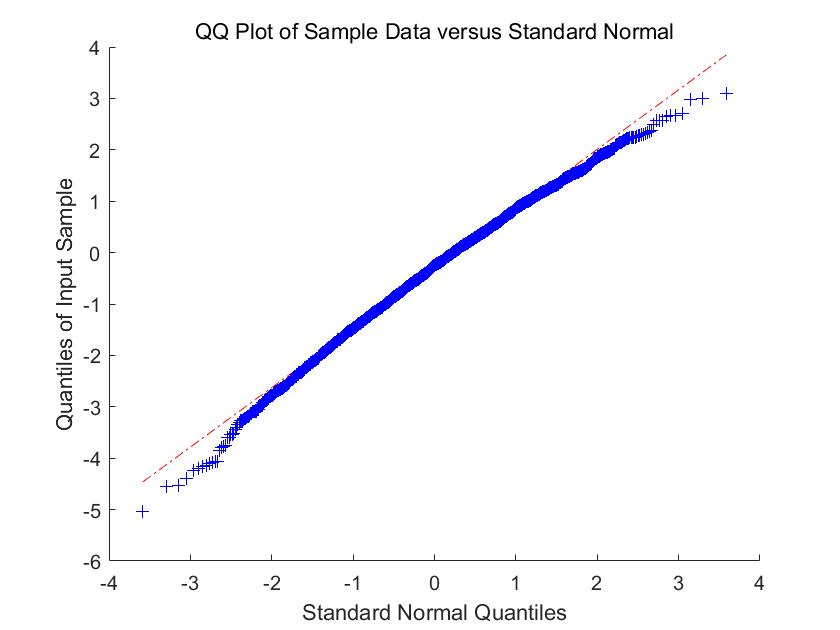}
	\caption{{The asymptotic normality of the second largest eigenvalue of noncentral Fisher matrix with $(p,n,N)=(200,2000,1000)$.}}
	\label{penG4}
\end{figure}

\begin{figure}[!h]	
	\centering
	\includegraphics[width=2.0in]{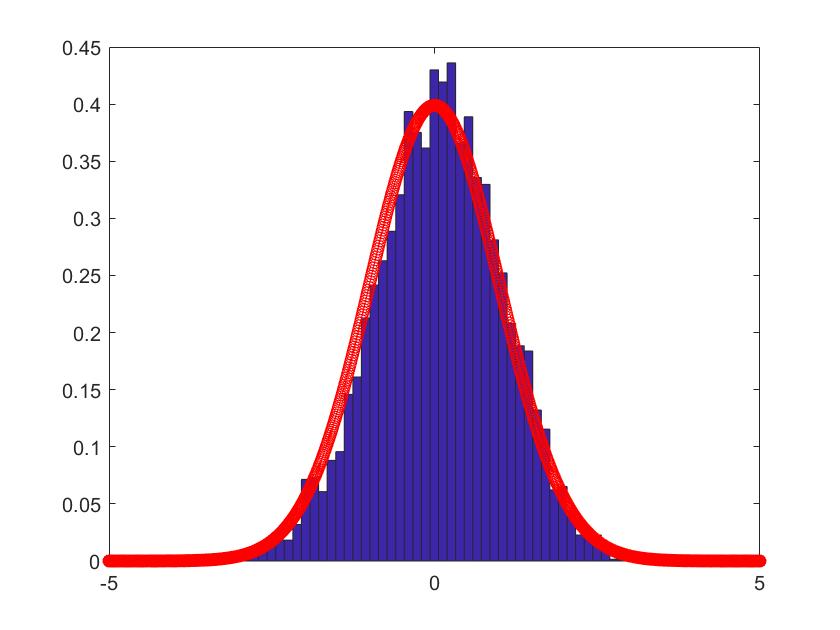}
	\includegraphics[width=2.0in]{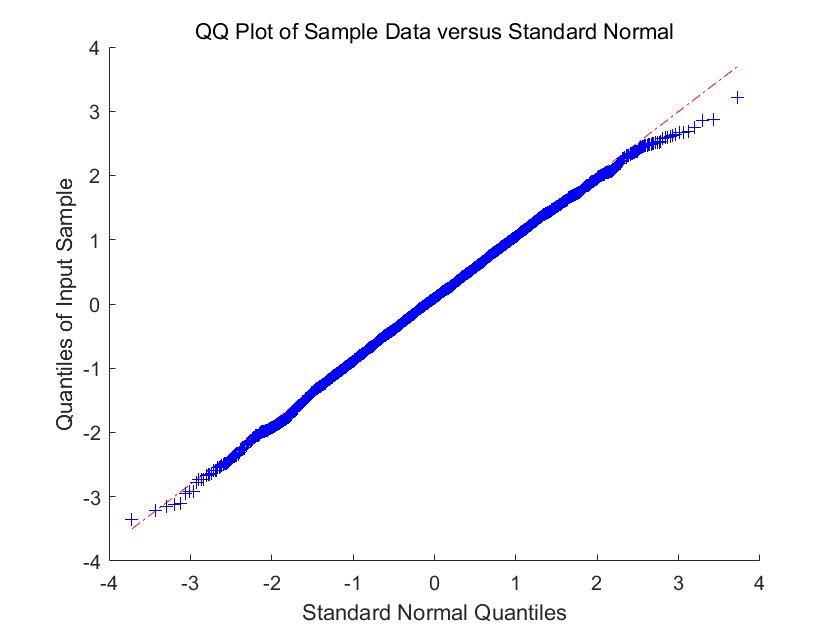}\\
	\caption{{The asymptotic normality of the largest eigenvalue of the CCA matrix with $(p,q,n)=(200,200,1000)$.}}
	\label{penG5}
\end{figure}

\begin{figure}[!h]	
	\centering
	\includegraphics[width=2.0in]{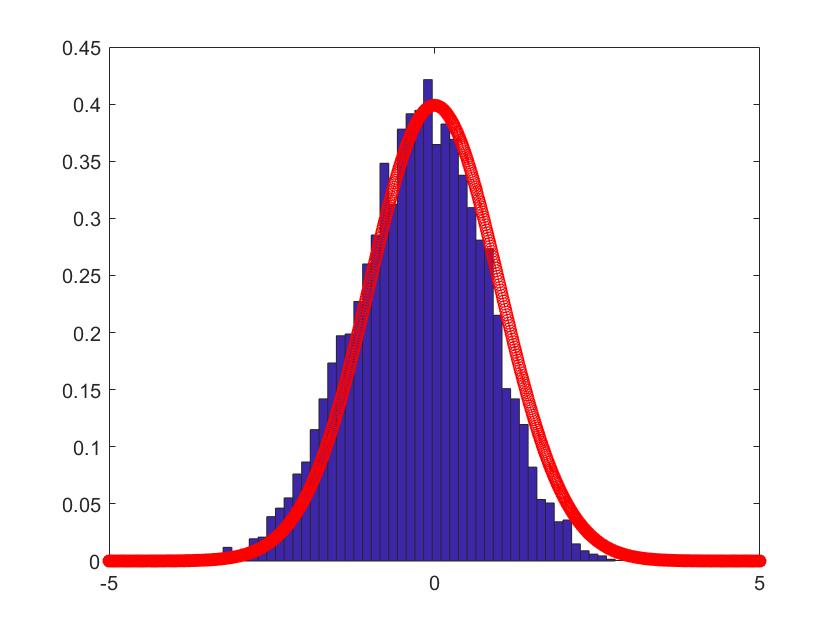}
	\includegraphics[width=2.0in]{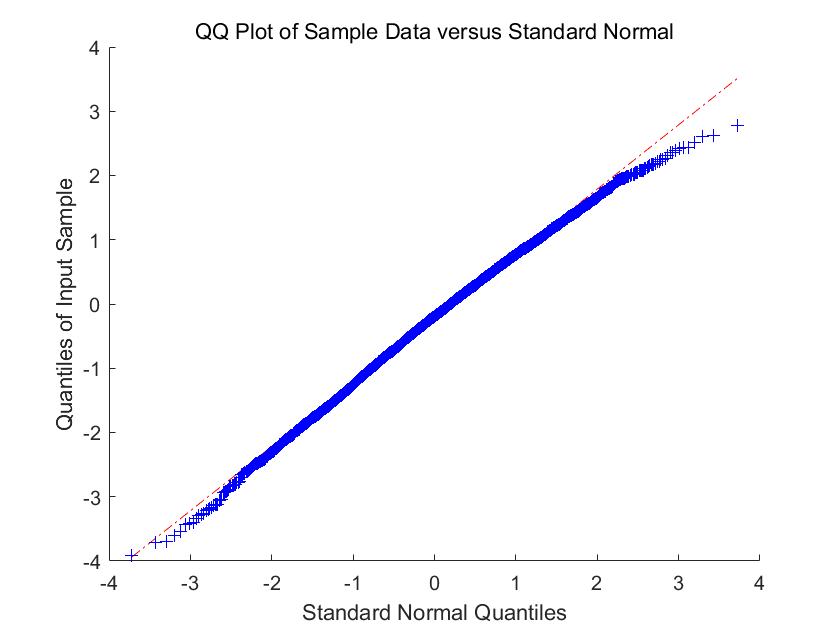}\\
	\caption{{The asymptotic normality of the second largest eigenvalue of the CCA matrix with $(p,q,n)=(200,200,1000)$.}}
	\label{penG6}
\end{figure}

In Figure \ref{penG1}-\ref{penG6}, we compare the empirical density (the blue histogram) of the two of largest eigenvalues of $\bbC_n$, $\bbF_p$ and the CCA matrix with the standard normal density curve (the red line) with 2000 repetitions. To be more convincing, we put the Q-Q plots together with the histograms. According to the histograms and Q-Q plots, we can conclude that the Theorem \ref{CLTtuS1}, \ref{CLTnonF} and  \ref{cltcca} are reasonable.
%
\subsection{Simulations for the estimators}
\begin{figure}[!h]
	\centering
	\includegraphics[width=0.325\textwidth]{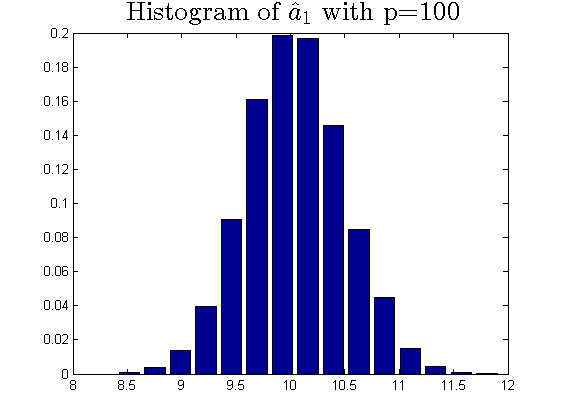}
	\includegraphics[width=0.325\textwidth]{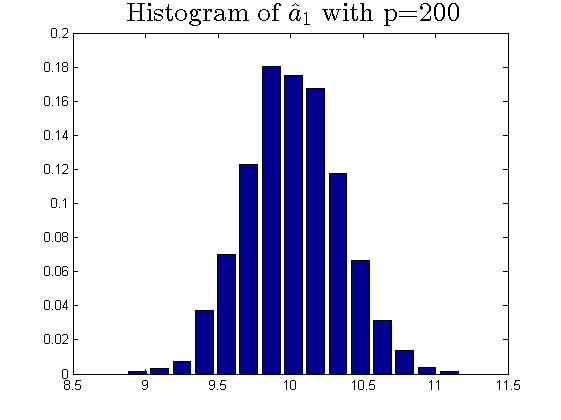}
	\includegraphics[width=0.325\textwidth]{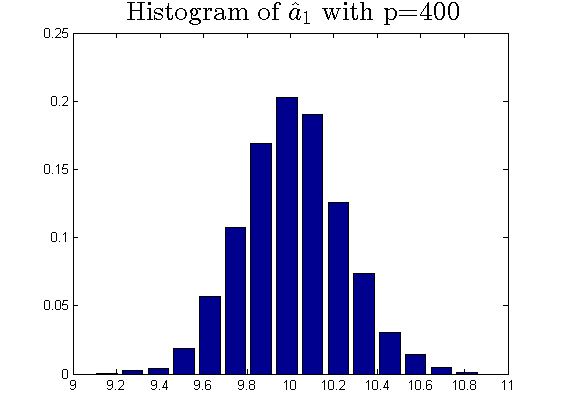}
	\caption{The estimator of $a_1$ $(a_1=10)$ by the results of  the noncentral sample covariance matrix with $p=100$, $200$, and $400$.}\label{f1}
\end{figure}

\begin{figure}[!h]
	\centering
	\includegraphics[width=0.325\textwidth]{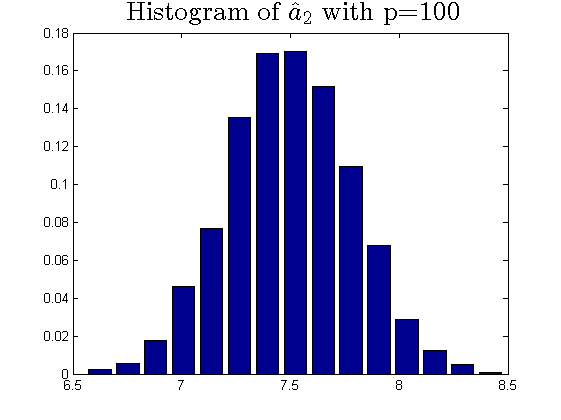}
	\includegraphics[width=0.325\textwidth]{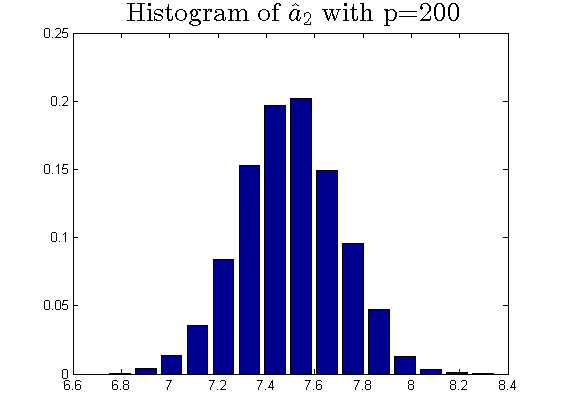}
	\includegraphics[width=0.325\textwidth]{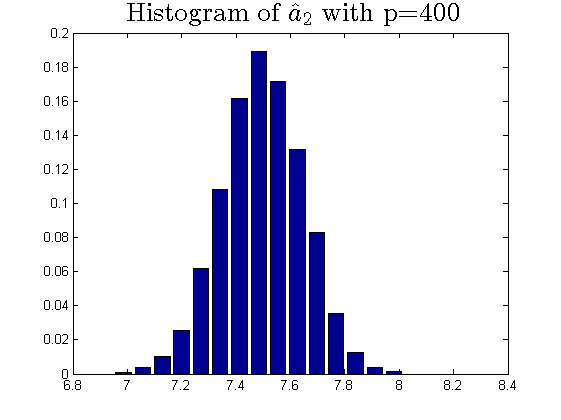}
	\caption{The estimator of $a_2$ $(a_2=7.5)$ by the results of the noncentral sample covariance matrix with $p=100$, $200$, and $400$.}\label{f2}
\end{figure}

\begin{figure}[!h]
	\centering
	\includegraphics[width=0.325\textwidth]{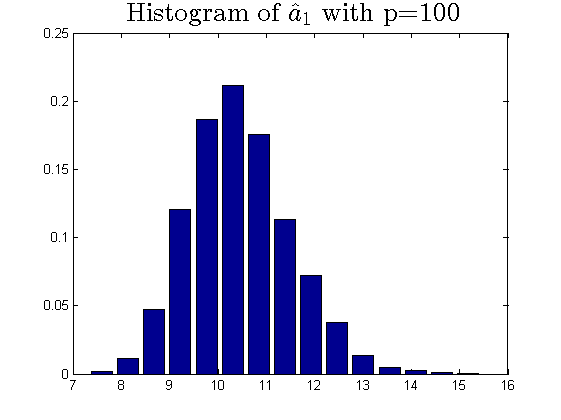}
	\includegraphics[width=0.325\textwidth]{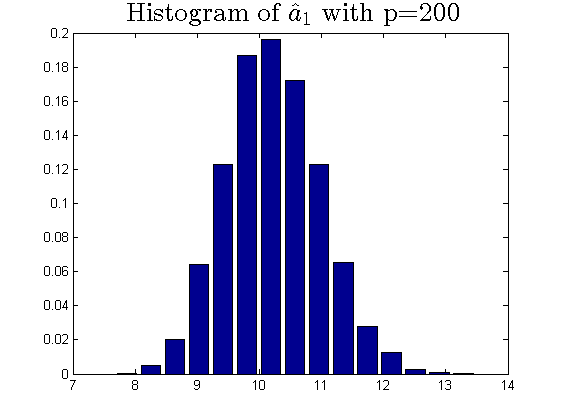}
	\includegraphics[width=0.325\textwidth]{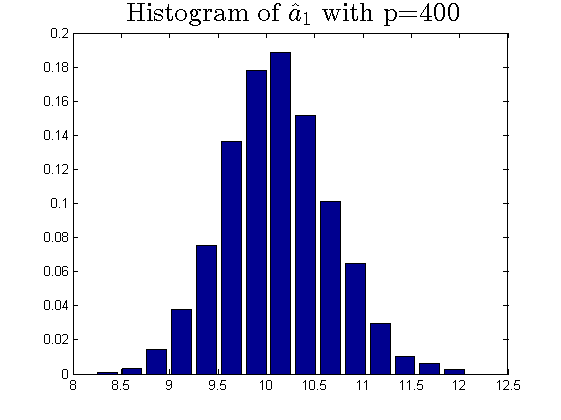}
	\caption{The estimator of $a_1$ $(a_1=10)$ by results of the noncentral Fisher matrix with $p=100$, $200$, and $400$.}\label{p1}
\end{figure}

\begin{figure}[!h]
	\centering
	\includegraphics[width=0.325\textwidth]{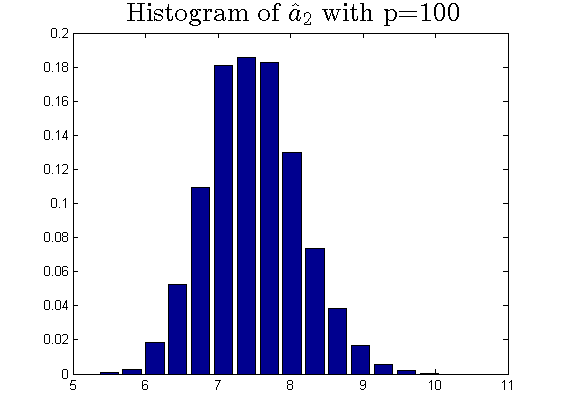}
	\includegraphics[width=0.325\textwidth]{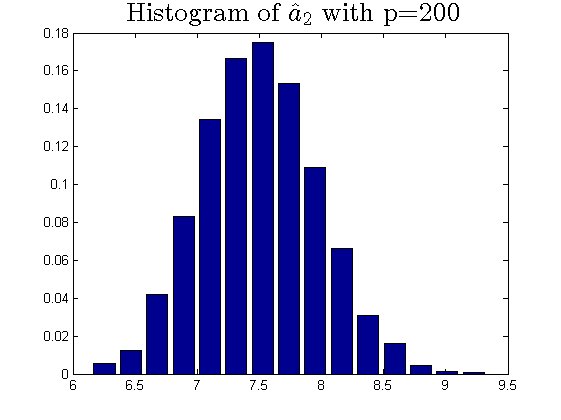}
	\includegraphics[width=0.325\textwidth]{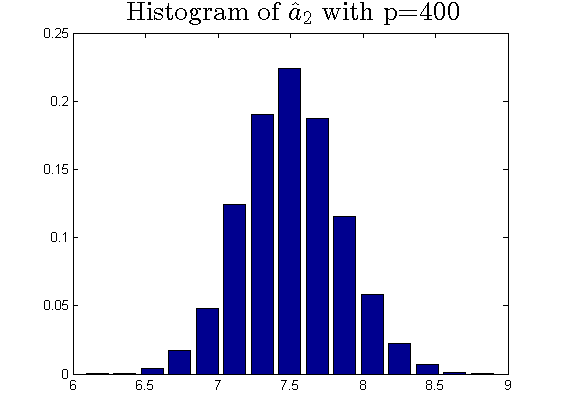}
	\caption{The estimator of  $a_2$ $(a_2=7.5)$ by results of the  noncentral Fisher matrix with $p=100$, $200$, and $400$.}\label{p2}
\end{figure}

\begin{figure}[!h]
	\centering
	\includegraphics[width=0.325\textwidth]{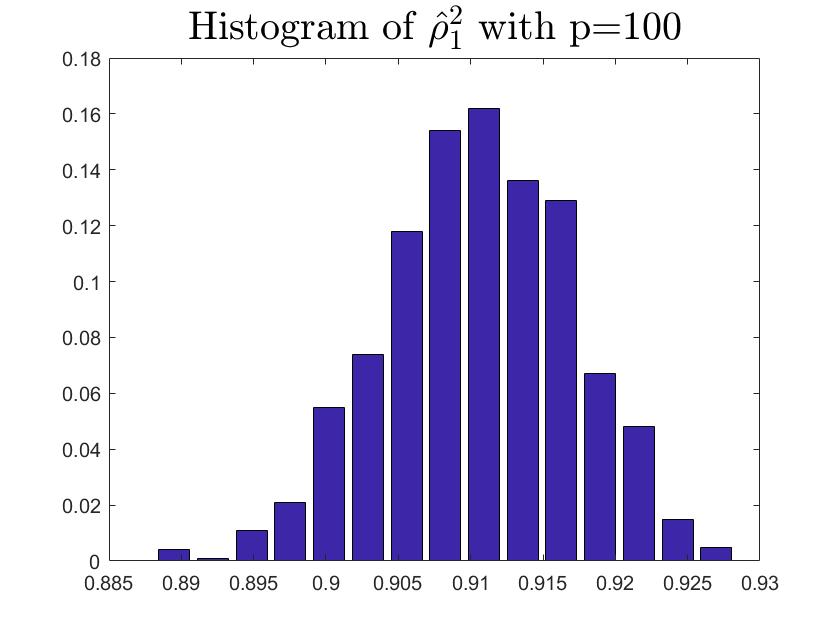}
	\includegraphics[width=0.325\textwidth]{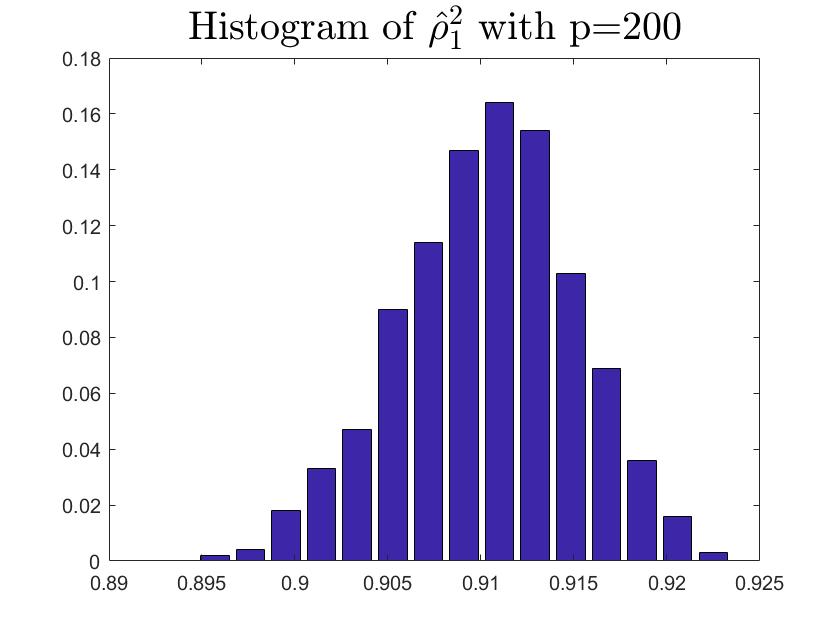}
	\includegraphics[width=0.325\textwidth]{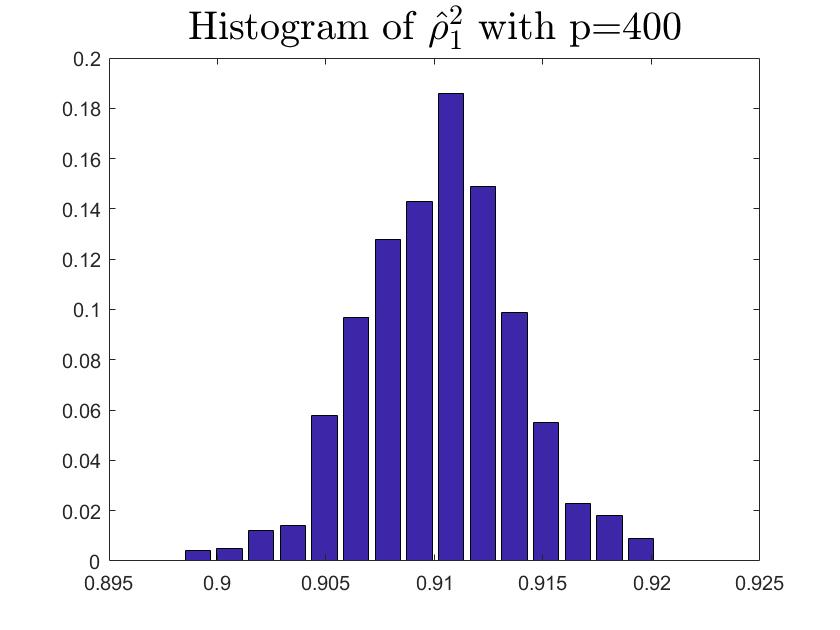}
	\caption{The estimator of $\rho_1^2$ $(\rho_1^2=10/11\approx0.9091)$ by results of the CCA matrix with $p=100$, $200$, and $400$.}\label{p3}
\end{figure}

\begin{figure}[!h]
	\centering
	\includegraphics[width=0.325\textwidth]{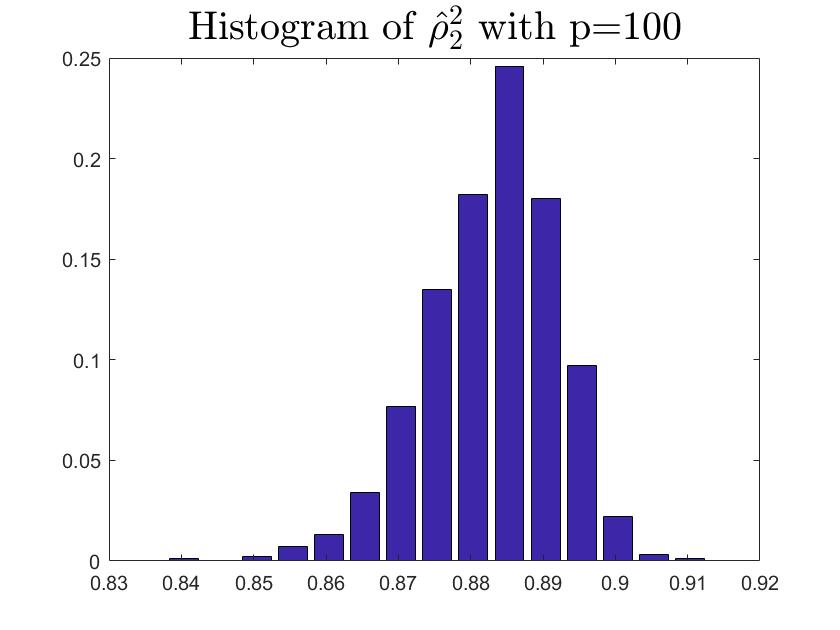}
	\includegraphics[width=0.325\textwidth]{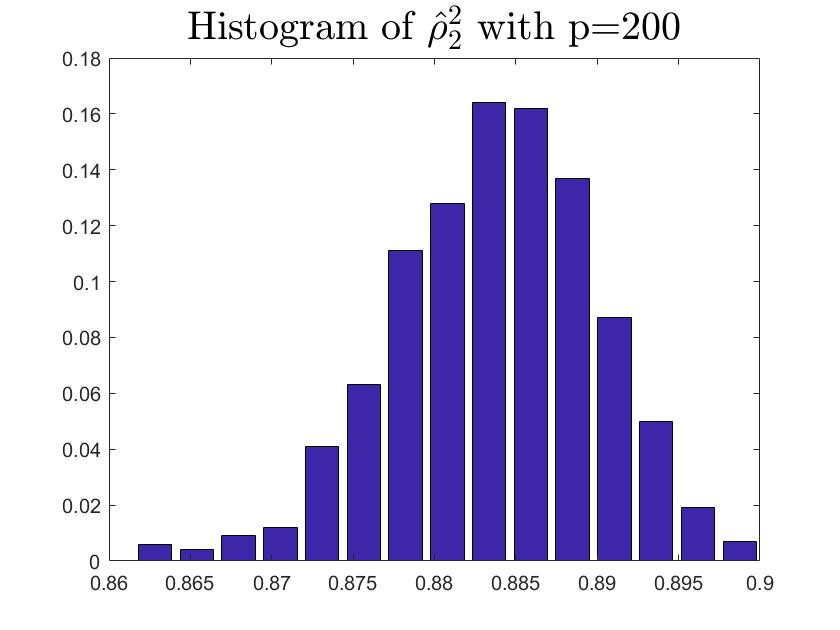}
	\includegraphics[width=0.325\textwidth]{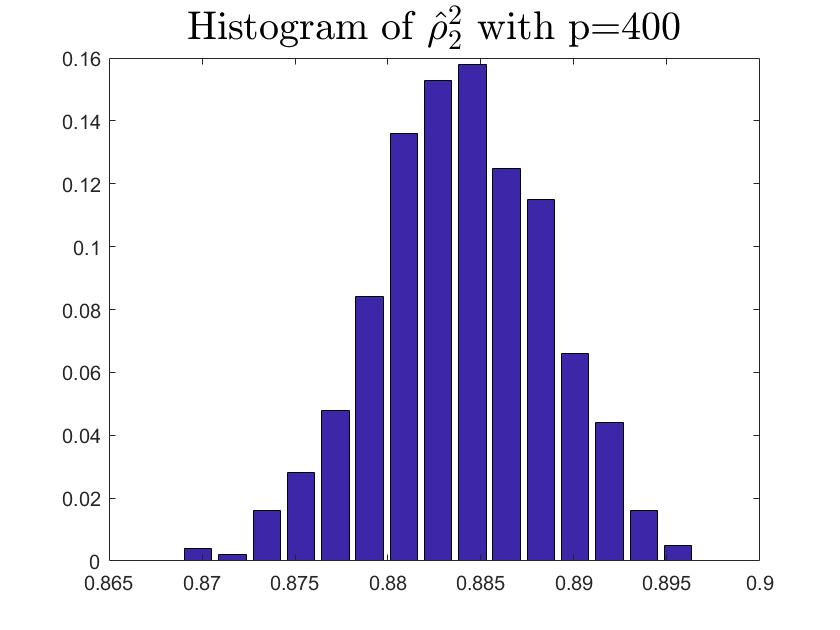}
	\caption{The estimator of $\rho_2^2$ $(\rho_2^2=15/17\approx0.8824)$ by results of the CCA matrix with $p=100$, $200$, and $400$.}\label{p4}
\end{figure}

We conduct  the following simulations to verify the accuracy of the estimators in (\ref{estimatorfornonS}), (\ref{estimatorfornonF}) and (\ref{estimatorforCCA}).  Unlike (\ref{Sim1model}), in this subsection, the eigenvalues of ${\mathbf \Xi}{\mathbf \Xi}^*/n$ are set as
\begin{align}\label{Simulationeigenvalues}
10, 7.5, 7.5, \underbrace{1,\cdots, 1}_{p-3},
\end{align}	
where $a_1=10$, $a_2=a_3=7.5$, and $H=\delta_{\{1\}}$.  Note that we set the second eigenvalue as a multiple eigenvalue. {According to  the single root condition in Assumption d$'$, we set the eigenvalues of $\boldsymbol{\Lambda}$ satisfying the model (\ref{Sim11model})}.
\begin{remark}
	{Compared with the assumption for the eigenvalues in (\ref{Simulationeigenvalues}), (\ref{Sim1model}) is set without multiple population spiked eigenvalues, which is considered that the joint distribution of the eigenvalues of multiple dimensional GOE matrix can not be visually displayed.}
\end{remark}
 We consider the estimator $\hat{a}_1$ and $\hat{a}_2$ with $p=100, 200$ and $400$, respectively. Then the frequency histograms of the estimators are present in Figure \ref{f1}-\ref{p2} with $5000$ repetitions. In Figure \ref{p3}-\ref{p4}, we show the accuracy of the estimator $\hat{\rho}_i$ for the two of largest population canonical correlation coefficients. We conclude that the estimators become accurate with $p$ increasing as the horizontal axis is more concentrated among three estimators. 

To a certain extent, the histograms present the accuracy of estimates of the spiked eigenvalues, however, which is not convincing enough. So we use the Mean Square Errors (MSE) criteria to show the accuracy of estimates in Table \ref{MSE}.
\begin{table}[!h]
	\newcommand{\tabincell}[2]{\begin{tabular}{@{}#1@{}}#2\end{tabular}}
	\renewcommand{\tabcolsep}{0.5pc}
	\renewcommand{\arraystretch}{1.5}
	\centering
	\caption{The MSE of the three estimators}
	\begin{tabular}{|c|ccc|ccc|}\hline
		&\multicolumn{3}{|c|}{$a_1=10$ {$\left(\rho_1=\sqrt{10/11}\right)$}}&\multicolumn{3}{|c|}{ $a_2=7.5$ {$\left(\rho_2=\sqrt{15/17}\right)$}}    \\ \hline
		$p$ & 100& 200 &400 &100& 200 &400\\ \hline
		$S$&1.2928 &0.6183 & 0.3176 &0.4017 &0.2194 & 0.1144\\ \hline
		$F$&0.2077&0.1064 & 0.0540 & 0.0821 & 0.0421& 0.0222 \\ \hline	
		{$CCA$}&4.8001e-05&2.5296e-05 &1.2784e-05&7.7782e-05 &4.4870e-05 &2.3276e-05 \\ \hline	
	\end{tabular}\label{MSE}
\end{table} 

The notation $S$, $F$ and $CCA$ in Table \ref{MSE} are short for the noncentral sample covariance matrix, the noncentral Fisher matrix and CCA matrix, respectively.  According to Table \ref{MSE}, we find that the MSE decreases as the dimensionality $p$ increases, which is consistent with the above mentioned simulation results. Compared with the MSE of the two estimators in Table \ref{MSE}, the estimator derived by the noncentral Fisher matrix is more consistent than the noncentral sample covariance matrix.

{\subsection{Simulation for multiple roots case}
Both Theorem (\ref{limitsCCA}) and Theorem (\ref{cltcca}) need satisfy the Assumption d$'$, specially, the spiked eigenvalue is a single root. So we set $\Lambda$ meet (\ref{Sim11model}) in subsection \ref{Snormal}. However, we guess our theories and estimator for the single root are effective under the multiple roots condition. In this subsection,  we will present some simulations to show the 
correctness and rationality estimating population canonical correlation coefficient under the multiple roots condition. We assume 
\begin{align}
\boldsymbol{\Lambda}=(\sqrt{10/11},\sqrt{15/17},\sqrt{15/17},\sqrt{1/2},\cdots,\sqrt{1/2}),
\end{align}
set the ratio among $(p:n:N)$ as $(1:3:9)$ and 1000 repetitions. According to Fig \ref{p5} and Table \ref{MSE2}, the performance looks good.
\begin{figure}[!h]
	\centering
	\includegraphics[width=0.325\textwidth]{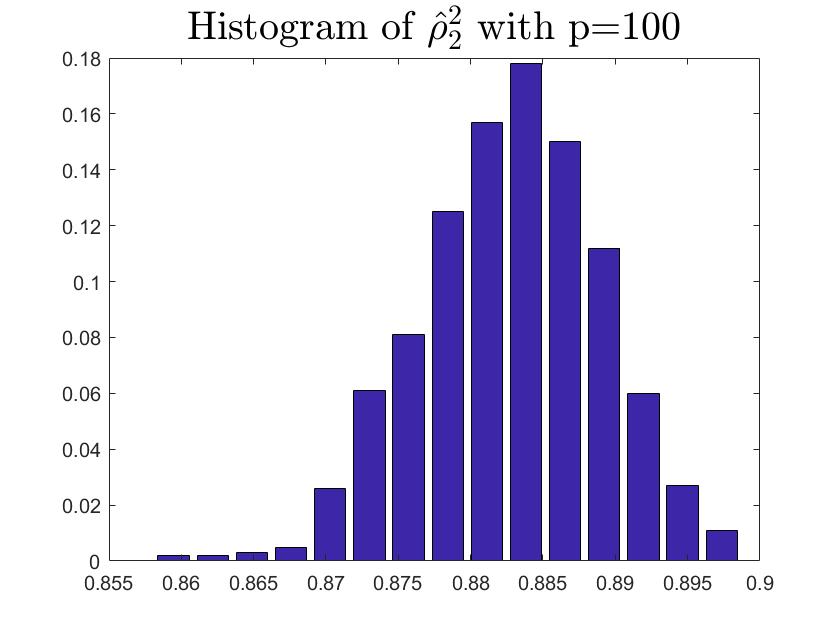}
	\includegraphics[width=0.325\textwidth]{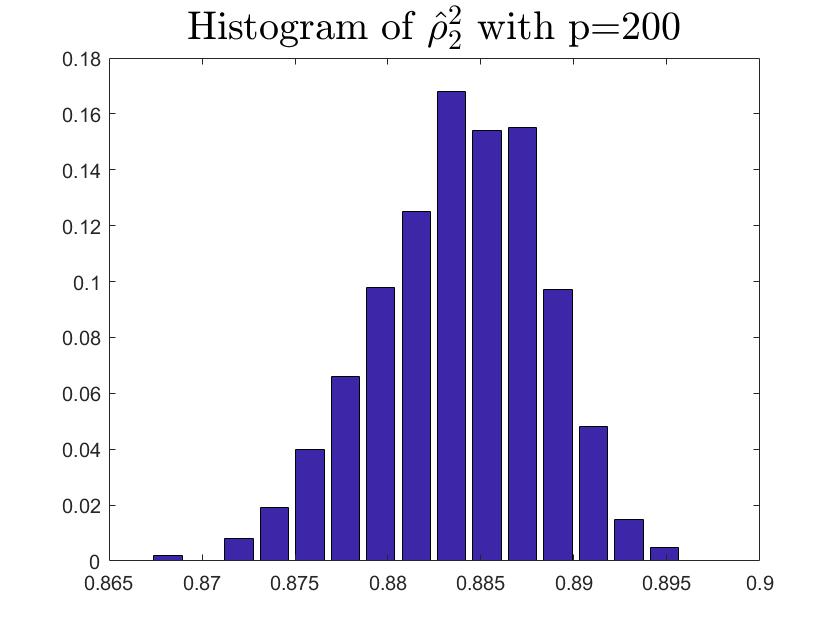}
	\includegraphics[width=0.325\textwidth]{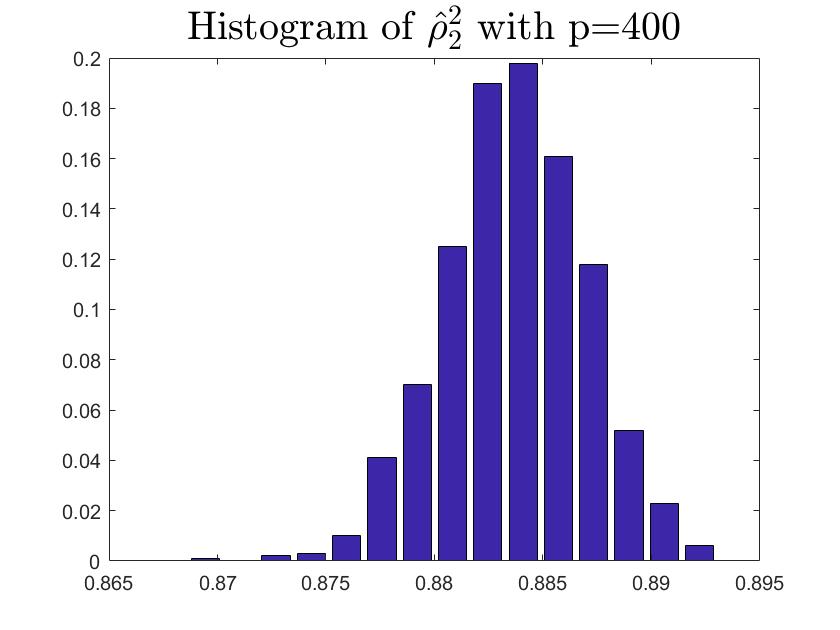}
	\caption{The estimator of $\rho_2^2$ $(\rho_2^2=15/17\approx0.8824)$ by results of the CCA matrix with $p=100$, $200$, and $400$.}\label{p5}
\end{figure}
\begin{table}[!h]
	\newcommand{\tabincell}[2]{\begin{tabular}{@{}#1@{}}#2\end{tabular}}
	\renewcommand{\tabcolsep}{0.8pc}
	\renewcommand{\arraystretch}{1.5}
	\centering
	\caption{The MSE of the population ccc estimators under multiple roots}
	\begin{tabular}{|c|ccc|}\hline
		&\multicolumn{3}{|c|}{ {$\rho_2=\rho_3=\sqrt{15/17}$}}   \\ \hline
		$p$ & 100& 200 &400\\ \hline	
		{$CCA$}&3.7055e-05&2.2209e-05 &1.2577e-05\\ \hline	
	\end{tabular}\label{MSE2}
\end{table} 
}

%

\section{Conclusion and discussion}
In this work, we study the limiting properties of  the sample spiked eigenvalues of the noncentral Fisher matrix under high-dimensional setting and Gaussian assumptions. Similar to the spiked model work, we find a phase transition for the limiting properties of the sample spiked eigenvalues of the noncentral Fisher matrix. In addition, we present the CLT for the sample spiked eigenvalues. As an accessory to the proof of the results of the limiting properties of the sample spiked eigenvlaues of the noncentral Fisher matrix, the fluctuations of the spiked eigenvalues of noncentral sample covariance matrix ${\mathbf C}_n$ are studied, which should have its own interests.

{\it General distribution.} It would be natural to ask if our results in the current work can be extended from Gaussian to more general distribution of the matrix entries. Our future work will focus on the limiting properties of the sample spiked eigenvalues of the noncentral Fisher matrix under general distribution of the matrix entries.

\section{Appendix}
\subsection{ Proof of Theorem \ref{limittuS1}}
There is no loss of generality in assuming  
\begin{eqnarray}\label{eigen-decom}
	\boldsymbol\Xi=\left(\begin{array}{c}\boldsymbol\Xi_1\\ \boldsymbol\Xi_2\end{array}\right)=\left(\begin{array}{cc} \bbD_{11} &\textbf{0}\\
		\textbf{0}& \tilde{\bbD}_{22} \end{array} \right)=\left(\begin{array}{ccc} \bbD_{11} &\textbf{0} &\textbf{0} \\
		\textbf{0}& \bbD_{22} &\textbf{0} \end{array} \right),
\end{eqnarray}
where $\boldsymbol\Xi_1$ is the first $M$ rows of $\boldsymbol\Xi$,
the diagonal of $\bbD_{11}$ is composed of the $M$ spiked eigenvalues
while the diagonal of $\bbD_{22}$ is composed of the non-spiked bounded eigenvalues, and $\tilde{\bbD}_{22}=(\bbD_{22}, \textbf{0})$.
According to the structure of $\boldsymbol\Xi$ in (\ref{eigen-decom}), we decompose $\bbX$ as following,
\begin{eqnarray*}
	\bbX=\left(\begin{array}{c} \bbX_1\\ \bbX_2 \end{array} \right),
\end{eqnarray*}
where $\bbX_1$ denotes the first $M$ rows of $\bbX$, and $\bbX_2$ stands for the remaining rows of $\bbX$.
Then the sample eigenvalues of noncentral sample  covariance matrix $\bbC_n$ are sorted in descending order as
\begin{eqnarray*}
	l_1^{\bbC_n} \geq l_2^{\bbC_n}\geq\cdots \geq l_p^{\bbC_n}.
\end{eqnarray*}
If we only consider the sample spiked eigenvalues of $\bbC_n$, $l_i^{\bbC_n}$, $i\in\mathcal{J}_k$, $k=1,\cdots,K$, then the eigen-equation is
\begin{align}\label{che}
	0=&\left|l_i^{\bbC_n}\bbI_p-\frac{1}{n}(\boldsymbol\Xi+\bbX)(\boldsymbol\Xi+\bbX)^{\ast}\right|\nonumber\\
	=&\left|l_i^{\bbC_n}\bbI_p-\frac{1}{n}\left(\begin{array}{c}\boldsymbol\Xi_1+\bbX_1 \\
		\boldsymbol\Xi_2+\bbX_2 \end{array}\right)\left((\boldsymbol\Xi_1+\bbX_1)^{\ast},(\boldsymbol\Xi_2+\bbX_2)^{\ast}\right)\right|\\
   =&\left|\begin{array}{cc}l_i^{\bbC_n}\bbI_M\!-\!\frac{1}{n}(\boldsymbol\Xi_1\!+\!\bbX_1)(\boldsymbol\Xi_1\!+\!\bbX_1)^{\ast} & \!-\frac{1}{n}(\boldsymbol\Xi_1\!+\!\bbX_1)(\boldsymbol\Xi_2\!+\!\bbX_2)^{\ast}\\
		-\frac{1}{n}(\boldsymbol\Xi_2\!+\!\bbX_2)(\boldsymbol\Xi_1\!+\!\bbX_1)^{\ast} &
		l_i^{\bbC_n}\bbI_{p-M}-\frac{1}{n}(\boldsymbol\Xi_2\!+\!\bbX_2)(\boldsymbol\Xi_2\!+\!\bbX_2)^{\ast}\end{array}
	\right|.\nonumber
\end{align}
Because $M$ is fixed,  $\boldsymbol\Xi_2\boldsymbol\Xi_2^*/n$ and  $\boldsymbol\Xi\boldsymbol\Xi^*/n$ share the same the LSD, then $l_i^{\bbC_n}$ is an outlier for large $n$, i.e. $| l_i^{\bbC_n} \bbI_{p-M}-\frac{1}{n}(\boldsymbol\Xi_2+\bbX_2 )(\boldsymbol\Xi_2^*+\bbX_2^*) |\neq  0$. Rewrite (\ref{che}) by the inverse of a partitioned matrix (\cite{HJ2012M}, Section 0.7.3) and the in-out exchange formula, we have
\begin{align*}
	0=&|l_i^{\bbC_n}\bbI_M-\frac{1}{n}(\boldsymbol\Xi_1+\bbX_1)(\boldsymbol\Xi_1+\bbX_1)^{\ast}-\frac{1}{n}(\boldsymbol\Xi_1+\bbX_1)(\boldsymbol\Xi_2+\bbX_2)^{\ast}\\
	&\times\left[l_i^{\bbC_n}\bbI_{p-M}-\frac{1}{n}(\boldsymbol\Xi_2+\bbX_2)(\boldsymbol\Xi_2+\bbX_2)^{\ast}\right]^{-1}\frac{1}{n}(\boldsymbol\Xi_2+\bbX_2)(\boldsymbol\Xi_1+\bbX_1)^{\ast}|\\
	\Longleftrightarrow0\!=\!&\left|\bbI_M\!-\!\frac{1}{n}(\boldsymbol\Xi_1+\bbX_1)\left(l_i^{\bbC_n}\bbI_n\!-\!\frac{1}{n}(\boldsymbol\Xi_2\!+\!\bbX_2)^{\ast}(\boldsymbol\Xi_2\!+\!\bbX_2)\right)^{-1}(\boldsymbol\Xi_1\!+\!\bbX_1)^{\ast}\right|,
\end{align*}
if $l_i^{\bbC_n}\neq 0$.
For simplicity, we denote $[\frac{1}{n}(\boldsymbol\Xi_2\!+\!\bbX_2)^{\ast}(\boldsymbol\Xi_2\!+\!\bbX_2)-l_i^{\bbC_n}\bbI_n]^{-1}$ briefly by $\bbA_n(l_i^{\bbC_n})$. When there is no ambiguity, we also rewrite it as $\bbA$.
Then
\begin{align}
	\boldsymbol{\Omega}_n^{\bbC_n}\!\overset{\triangle}{=}&\bbI _M + \frac{1}{n}(\boldsymbol\Xi_1 + \bbX_1 )\bbA_n(l_i^{\bbC_n})(\boldsymbol\Xi_1^* + \bbX_1^* ) \nonumber\\
	=&\bbI_M +\frac1n({\rm tr}\bbA_n(l_i^{\bbC_n}))\bbI_M -\frac{1}{l_i^{\bbC_n}(1+c_{1n}m_{2n}(l_i^{\bbC_n}))}\frac{1}{n}\boldsymbol\Xi_1\boldsymbol\Xi_1^* + \boldsymbol{\Omega}_0^{\bbC_n},
	\label{eqom}
\end{align}
where
\begin{align}\label{305}
	{\boldsymbol\Omega}_0^{\bbC_n}(l_i^{\bbC_n})\!&=\!\frac1n\bbX_1\bbA\bbX_1^*\!-\!\frac1n({\rm tr}\bbA)\bbI_M\!+\!\frac1n(\bbX_1\bbA\boldsymbol\Xi_1^* 
	\!+\!\boldsymbol\Xi_1\bbA\bbX_1^*)\nonumber\\
	\!&+\!\frac{1}{n}\bbD_{11}^*\left[\bbA_{11}\!+\!\frac{\bbI}{l_i^{\bbC_n}(1+c_{1n}m_{2n}(l_i^{\bbC_n}))}\right]\bbD_{11},\\
	m_{2n}(l_i^{\bbC_n})&=\frac1{p-M}{\rm tr}\left[\frac{1}{n}(\bbX_{22}\!+\!\tilde{\bbD}_{22})(\bbX_{22}\!+\!\tilde{\bbD}_{22})^*\!-\!l_i^{\bbC_n}\bbI_{p\!-\!M}\right]^{-1}
\end{align}
$c_{1n}\!=\!({p\!-\!M})/({n\!-\!M})$, and $\bbA_{11}$ is the first $M\times M$ major diagonal submatrix of $\bbA_n(l_i^{\bbC_n})$.

Here we announce that the following results of almost sure convergence are ture without proofs and the detailed proofs are postponed in Subsection \ref{523}.
\begin{eqnarray}\label{303}
&\bbA_{11}(l_i^{\bbC_n})+\frac{1}{l_i^{\bbC_n}+l_i^{\bbC_n}c_{1n}m_{2n}(l_i^{\bbC_n})}\bbI_M\overset{a.s.}{\longrightarrow}\boldsymbol{0}_{M\times M},\label{Slim1}\\
&\frac{1}{n}\bbX_1\bbA_n(l_i^{\bbC_n})\bbX_1^{\ast}-\frac{1}{n}(\tr\bbA_n(l_i^{\bbC_n}))\bbI_M\overset{a.s.}{\longrightarrow}\boldsymbol{0}_{M\times M}, \label{Slim2}\\
&\frac{1}{n}\bbX_1\bbA_n(l_i^{\bbC_n})\boldsymbol\Xi_1^{\ast}+\frac{1}{n}\boldsymbol\Xi_1\bbA_n(l_i^{\bbC_n})\bbX_1^{\ast}\overset{a.s.}{\longrightarrow}\boldsymbol{0}_{M\times M}\label{Slim3}.
\end{eqnarray}
{According to (\ref{eqom}), we {have} $\boldsymbol{\Omega}_0^{\bbC_n}\overset{a.s.}{\longrightarrow} \boldsymbol{0}_{M\times M}$.}
By the formula (1.3) in \cite{DS2007a}, we arrive at 
\begin{equation*}
\frac{1}{n}\tr\bbA_n(z)\overset{a.s.}{\longrightarrow}\underline{m}_2(z)\quad \mbox{uniformly for $z$}.
\end{equation*}
From
\begin{eqnarray}\label{eigeneq}
0=\left|\bbI_M+\frac{1}{n}\tr\bbA_n(l_i^{\bbC_n})\bbI-\frac{\frac{1}{n}\bbD_{11}\bbD_{11}^{\ast}}{l_i^{\bbC_n}+l_i^{\bbC_n}c_{1n}m_{2n}(l_i^{\bbC_n})}+\boldsymbol{\Omega}_0^{\bbC_n}(l_i^{\bbC_n})\right|,
\end{eqnarray}
we can obtain 
{\begin{align}\label{59}
	\left|\bbI_M+\underline{m}_{2n}^0(l_i^{\bbC_n})\bbI_M-\frac{\frac{1}{n}\bbD_{11}^{-1}\bbD_{11}^{-1}}{l_i^{\bbC_n}(1+c_{1n}m_{2n}^0(l_i^{\bbC_n}))}\right|\overset{a.s.}{\longrightarrow} 0.
	\end{align}}
For arbitrary $k$, let 
\begin{align}\label{lkcdef}
\lambda_{nk}^{\bbC}\overset{\triangle}{=}a_{k}\left(1-c_{1n}m_{1n}^0(a_{k})\right)^2+(1-c_{1n})(1-c_{1n}m_{1n}^0(a_{k}))
\end{align}
where $c_{1n}=(p-M)/(n-M)$, and $m_{1n}^0$ denotes the ST of ESD of $p-M$ bulk eigenvalues of $\boldsymbol{\Xi}\boldsymbol{\Xi}^*/n$.
An easy calculation shows that 
\begin{align*}
1+\underline{m}_{2n}^0(\lambda_{nk}^{\bbC})-\frac{a_k}{\lambda^{\bbC}_k(1+c_1m_{2n}^0(\lambda^{\bbC}_{nk}))}=0,
\end{align*}
where $m_{2n}^0$ is the ST of LSD of $(\boldsymbol{\Xi}_2+\bbX_2)(\boldsymbol{\Xi}_2+\bbX_2)^*/n$ with $c_1$ and $H$ replaced by $p-M/(n-M)$ and ESD of $p-M$ bulk eigenvalues of $\boldsymbol{\Xi}\boldsymbol{\Xi}^*/n$.
{
	Combining (\ref{59}) and the fact that the dimension of matrix is finite, there exists `$j$' (assume $j\in\mathcal{J}_k$) such that the $j$-th diagonal element convergence almost surely to zero.
	For this `$k$' we have} 
\begin{align*}
&\Big|\bbI_M\!+\!\underline{m}_{2n}^0(l_i^{\bbC_n})\bbI_M\!-\!\frac{\frac{1}{n}\bbD_{11}^{-1}\bbD_{11}^{-1}}{l_i^{\bbC_n}(1\!+\!c_{1n}m_{2n}^0(l_i^{\bbC_n}))}\\
&\!-\!\left[	1\!+\!\underline{m}_{2n}^0(\lambda_{nk}^{\bbC})\!-\!\frac{a_k}{\lambda^{\bbC}_{nk}(1\!+\!c_1m_{2n}^0(\lambda^{\bbC}_{nk}))}\right]\bbI_M\Big|\overset{a.s.}{\longrightarrow} 0.
\end{align*}
{Then subtracting the $j$-th one from all the diagonal elements of the matrix in the above determinant, we find the difference has a lower bound except the $k$-th block containing $j$-th location, i.e. } 
\begin{align*}
\frac{a_s-a_k}{l_i^{\bbC_n}(1+c_{1n}m_{2n}^0(l_i^{\bbC_n}))}, \quad s\notin\mathcal{J}_k
\end{align*}
{is with the lower bound as a result of (\ref{sc}) in Assumption A.} So for the diagonal elements of $k$-th block, we have 
\begin{align*}
&\underline{m}_{2n}^0(l_i^{\bbC_n})-\frac{a_k}{l_i^{\bbC_n}(1+c_{1n}m_{2n}^0(l_i^{\bbC_n}))}-\underline{m}_{2n}^0(\lambda_{nk}^{\bbC})+\frac{a_k}{\lambda_{nk}^{\bbC}(1+c_{1n}m_{2n}^0(\lambda_{nk}^{\bbC}))}\overset{a.s.}{\longrightarrow} 0\\
\Leftrightarrow&\left(\frac{l_i^{\bbC_n}-\lambda_{nk}^{\bbC}}{\lambda_{nk}^{\bbC}}\right)\left[\lambda_{nk}^{\bbC}(\underline{m}_{2n}^0)'(\xi_1)+\frac{a_k}{l_i^{\bbC_n}}\frac{1+c_{1n}m_{2n}^0(\lambda_{nk}^{\bbC})+l_i^{\bbC_n}(m_{2n}^0)'(\xi_2)}{(1+c_{1n}m_{2n}^0(l_i^{\bbC_n}))(1+c_{1n}m_{2n}^0(\lambda_{nk}^{\bbC}))}\right]\overset{a.s.}{\longrightarrow} 0, 
\end{align*}
By the factor 
\begin{align*}
\lambda_{nk}^{\bbC}(\underline{m}_{2n}^0)'(\xi_1)+\frac{a_k}{l_i^{\bbC_n}}\frac{1+c_{1n}m_{2n}^0(\lambda_{nk}^{\bbC})+l_i^{\bbC_n}(m_{2n}^0)'(\xi_2)}{(1+c_{1n}m_{2n}^0(l_i^{\bbC_n}))(1+c_{1n}m_{2n}^0(\lambda_{nk}^{\bbC}))}
\end{align*}
{has lower bound},
we get $(l_i^{\bbC_n}-\lambda_{nk}^{\bbC})/\lambda_{nk}^{\bbC}\overset{a.s.}{\longrightarrow} 0$. {If $a_k$ is bounded, the limit $\lambda_k^{\bbC}$ of $\lambda_{nk}^{\bbC}$ satisfies}
\begin{eqnarray*}
	&&0=1-\frac{1-c_1}{\lambda^{\bbC}_k}+c_1m_2(\lambda^{\bbC}_k)+\frac{a_k}{-{\lambda^{\bbC}_k}-\lambda^{\bbC}_kc_1m_2(\lambda^{\bbC}_k)}\\
	\Longleftrightarrow&&a_k=\lambda^{\bbC}_k\left(1+c_1m_2(\lambda^{\bbC}_k)\right)^2-(1-c_1)\left(1+c_1m_2(\lambda^{\bbC}_k)\right)\\
	\Longleftrightarrow&&a_k=\left[\lambda^{\bbC}_k\left(1+c_1\frac{m_1(a_k)}{1-c_1m_1(a_k)}\right)-(1-c_1)\right]\left[1+c_1\frac{m_1(a_k)}{1-c_1m_1(a_k)}\right]\\
	\Longleftrightarrow&&\psi_{\mathbf C}(a_k)\overset{\bigtriangleup}{=}\lambda^{\bbC}_k=a_k\left(1-c_1m_1(a_k)\right)^2+(1-c_1)\left(1-c_1m_1(a_k)\right),
\end{eqnarray*}
where $m_1$ is the Stieltjes transform of the LSD $H$ of $\boldsymbol\Xi\boldsymbol\Xi^{\ast}/n$.
The second equivalence relation is a consequence of (\ref{m_2}). The proof of Theorem \ref{limittuS1} is completed.

\begin{remark}
	{If $a_k$ tends to infinite as stated in Assumption A, we only need to  multiply by $\sqrt{n\cdot l_i}\bbD_{11}^{-1}$ from both side of (\ref{eqom}), and similar arguments above applying to  the infinite case, we can get the same conclusion with only describing as: $(l_i^{\bbC_n}-\lambda_{nk}^{\bbC})/\lambda_{nk}^{\bbC}\overset{a.s.}{\longrightarrow} 0$. Therefore, we will not repeat the situation that $a_k$ tending infinite in following proof.}
\end{remark}

\subsubsection{Proofs of (\ref{Slim1}), (\ref{Slim2}) and (\ref{Slim3})}\label{523}
(\ref{Slim3}) can be obtain by Kolmogorov's law of large numbers straightly. The proof of {(\ref{Slim1}) and (\ref{Slim2}) are similar, so we take the proof of (\ref{Slim2}) as example.} We consider {the following series for any $\varepsilon>0$},
\begin{align*}
	&\sum_{n=1}^{\infty}\rP(\|\frac{1}{n}\bbX_1\bbA_n\bbX_1^*-\frac{1}{n}(\tr\bbA)\bbI_M\|_K>\varepsilon)\\
	=&\sum_{n=1}^{\infty}\rP(\|\frac{1}{n}\bbX_1\bbA_n\bbX_1^*-\frac{1}{n}(\tr\bbA)\bbI_M\|_K>\varepsilon,\mathcal{A})\\
	&+\sum_{n=1}^{\infty}\rP(\|\frac{1}{n}\bbX_1\bbA_n\bbX_1^*-\frac{1}{n}(\tr\bbA)\bbI_M\|_K>\varepsilon,\mathcal{A}^c)
\end{align*}
where the event $\mathcal{A}$ means the spectral norm of $\bbA$ is bounded, i.e., $\|\bbA\|\leq C$ and $\|\cdot\|_K$ means Kolmogorov norm, defined as the largest absolute value of all the entries. Then we have 
\begin{align*}
		\rP(\|\frac{1}{n}\bbX_1\bbA_n\bbX_1^*-\frac{1}{n}(\tr\bbA)\bbI_M\|_K>\varepsilon,\mathcal{A}^c)\leq \rP(\mathcal{A}^c)=o(n^{-t}),
	\end{align*}
	where the last equation is a consequence of the exact separation result of information-plus-noise type matrices in \cite{BS2012}.
For the first term, we have 
\begin{align}\label{summable}
	&\sum_{n=1}^{\infty}\rP(\|\frac{1}{n}\bbX_1\bbA_n\bbX_1^*-\frac{1}{n}(\tr\bbA)\bbI_M\|_K>\varepsilon\cap\mathcal{A})\nonumber\\
	=&\sum_{n=1}^{\infty}\E I\{\|\frac{1}{n}\bbX_1\bbA_n\bbX_1^*-\frac{1}{n}(\tr\bbA)\bbI_M\|_K>\varepsilon\cap\mathcal{A}\}\nonumber\\
	=&\sum_{n=1}^{\infty}\E[\E I\{\|\frac{1}{n}\bbX_1\bbA_n\bbX_1^*-\frac{1}{n}(\tr\bbA)\bbI_M\|_K>\varepsilon\cap\mathcal{A}\}|\bbA]\nonumber\\
	=&\sum_{n=1}^{\infty}\E[\rP(\|\frac{1}{n}\bbX_1\bbA_n\bbX_1^*-\frac{1}{n}(\tr\bbA)\bbI_M\|_K>\varepsilon\cap\mathcal{A})|\bbA]\nonumber\\
	\leq&\sum_{n=1}^{\infty}\E[\frac{1}{\varepsilon^{4r}}\E\|\frac{1}{n}\bbX_1\bbA_n\bbX_1^*-\frac{1}{n}(\tr\bbA)\bbI_M\|_K^{4r}I_{\mathcal{A}}|\bbA].
\end{align}

In fact, there exists a constant $K$, independent with $i,l,n$, {such that}  
\begin{align*}
	\E|(\frac{1}{n}\bbX_1\bbA_n\bbX_1^*-\frac{1}{n}(\tr\bbA)\bbI_M)_{il}|^{4}I_{\mathcal{A}}|\bbA\leq \frac{K}{n^2}
\end{align*}
which implies (\ref{summable}) is summable when $r=1$. By the Borel-Cantelli lemma, we have  
\begin{align*}
	\frac{1}{n}\bbX_1\bbA_n\bbX_1^*-\frac{1}{n}(\tr\bbA)\bbI_M\overset{a.s.}{\to}\boldsymbol{0}_{M\times M}.
\end{align*}

\subsection{Proof of Theorem \ref{CLTtuS1}}\label{CLT2}

In this section, we will consider the random vector 
\begin{align*}
	\gamma_k^{\bbC_n}=\sqrt{n}\{l_i^{\bbC_n}/\lambda_{nk}^{\bbC}-1,i\in\mathcal{J}_k\},
\end{align*}
where $\lambda_{nk}^{\bbC}$ is defined as (\ref{lkcdef}), 
The {reason} of using $\lambda_{nk}^{\bbC}$ rather than its limit $\lambda_{k}^{\bbC}$ lies in the fact that the convergence may be very slow.
The following proof is based on (\ref{eqom}) and (\ref{305}), then
\begin{equation*}
	0\!=\!\left|\bbI_M\!+\!\underline{m}_{2n}^0(\lambda_{nk}^{\bbC})\bbI\!-\!\frac{\frac{1}{n}\bbD_{11}\bbD_{11}^{\ast}}{\lambda_{nk}^{\bbC}(1\!+\!c_{1n}m_{2n}^0(\lambda_{nk}^{\bbC}))}\!+\!\boldsymbol{\Omega}_0^{\bbC_n}(\lambda_{nk}^{\bbC})\!+\!\varepsilon_1\bbI_M\!+\!\varepsilon_2\frac{1}{n}\bbD_{11}\bbD_{11}^*\!+\!\varepsilon_3\right|,
\end{equation*}
where $m_{2n}^0$ is the Stieltjes transform of $F^{\bbC}$ with parameter $H$ and $c_1$ replaced by $H_n$ and $c_{1n}$, and 
\begin{eqnarray}
	\label{e1}	\varepsilon_1&=&\frac{1}{n}\tr\bbA_n(l_i^{\bbC_n})-\underline{m}_{2n}^0(\lambda_{nk}^{\bbC})\\
	\label{e2}	\varepsilon_2&=&\frac{-1}{l_i^{\bbC_n}+l_i^{\bbC_n}\frac{p-M}{n}m_{2n}(l_i^{\bbC_n})}-\frac{-1}{\lambda_{nk}^{\bbC}+\lambda_{nk}^{\bbC}c_{1n}m_{2n}^0(\lambda_{nk}^{\bbC})},\\
	\label{e3}	\varepsilon_3&=&\boldsymbol{\Omega}_0^{\bbC_n}(l_i^{\bbC_n})-\boldsymbol{\Omega}_0^{\bbC_n}(\lambda_{nk}^{\bbC}).
\end{eqnarray}
Here we give the estimators of $\varepsilon_1$, $\varepsilon_2$, and $\varepsilon_3$ and the detailed proof is postponed in the following part,
\begin{eqnarray}
	\label{401}	\varepsilon_1&=&\frac{\gamma_{k}^{\bbC_n}}{\sqrt{n}}\lambda_k^{\bbC_n}[\underline m_2'(\lambda_k^{\bbC_n})+o_p(1)],\\
	\varepsilon_2&=&\frac{\gamma_{k}^{\bbC_n}}{\sqrt{n}}\left[\frac{1+c_1m_2(\lambda_k^{\bbC_n})+c_1\lambda_k^{\bbC_n}m_2'(\lambda_k^{{\bbC_n}})}{\lambda_k^{\bbC_n}[1+c_1m_2(\lambda_k^{\bbC_n})]^2}+o_p(1)\right],\\
	\varepsilon_3&=&o_p\left(\frac{1}{\sqrt{n}}\right)\boldsymbol{1}\boldsymbol{1}'.
\end{eqnarray}

According to the definition (\ref{lkcdef}), 
if $i\in\mathcal{J}_k$, then we obtain 
\begin{equation}
	1+\underline{m}_{2n}^0(\lambda_{nk}^{\bbC})-\frac{a_k}{\lambda_{nk}^{\bbC}(1+c_{1n}m_{2n}^0(\lambda_{nk}^{\bbC}))}=0.
\end{equation}
We rewrite the $k$-th block matrix of $\boldsymbol{\Omega}_n^{\bbC_n}$ as
\bqn
[\boldsymbol{\Omega}_n^{\bbC_n}]_{kk}=[\boldsymbol{\Omega}_0^{\bbC_n}(\lambda_{nk}^{\bbC})]_{kk}+\varepsilon_1\bbI_{m_k}+\varepsilon_2a_k\bbI_{m_k}+o_p(\frac{1}{\sqrt{n}}).
\eqn
By the discussion of the limiting distribution of $\boldsymbol{\Omega}_0^{\bbC_n}(\lambda_{nk}^{\bbC})$ and Skorokhod strong representation theorem, for more details, see \cite{Sk1956} or \cite{HB2014}, on an appropriate probability space, one may redefine the random variables such that $\boldsymbol{\Omega}_0^{\bbC_n}$ tends to the normal variables with probability one. 
Then, {the eigen-equation of (\ref{eqom})} becomes
\begin{align*}
0\!=\!\begin{vmatrix}	\frac{a_k(1\!-\!\frac{a_1}{a_k})}{\lambda_{nk}^{\bbC}b(\lambda_{nk}^{\bbC})}\!+\!O(n^{-1/2})&O(n^{-1/2})&O(n^{-1/2})\cr
	O(n^{-1/2})&\cdots&O(n^{-1/2})\cr
	O(n^{-1/2})&[\boldsymbol{\Omega}_0^{\bbC_n}]_{kk}\!+\!\varepsilon_1\bbI_{m_k}\!+\!\varepsilon_2a_k\bbI_{m_k}&\cdots\cr
	O(n^{-1/2})&\cdots&O(n^{-1/2})\cr
	O(n^{-1/2})&\cdots&\frac{a_k(1\!-\!\frac{a_M}{a_k})}{\lambda_{nk}^{\bbC}b(\lambda_{nk}^{\bbC})}\!+\!O(n^{-1/2})\cr\end{vmatrix}
\end{align*}
where $b(\lambda_{nk}^{\bbC})=1+c_{1n}m_{2n}^0(\lambda_{nk}^{\bbC})$ and $[\boldsymbol{\Omega}_0^{\bbC_n}]_{kk}$ is $k$-th diagonal block of $\boldsymbol{\Omega}_0^{\bbC_n}$.
Then multiplying $n^{1/4}$ to the $k$-th block row and column of the determinant of the eigen-equation above, and making $n\to\infty$. Then we have 
\begin{equation*}
	\sqrt{n}\left([\boldsymbol{\Omega}_0^{\bbC_n}]_{kk}+\varepsilon_1\bbI_{m_k}+\varepsilon_2a_k\bbI_{m_k}\right)\overset{{a.s.}}{\longrightarrow}0.
\end{equation*}
Simplifying the above, we have the random vector $\gamma_k^{\bbC_n}$ tends to a random vector consists of the ordered eigenvalues of GOE (GUE) matrix under real (complex) case with the scale parameter 
\begin{align}\label{theta1}
	\theta_1=&\frac{1}{[\lambda_k^{\bbC}\underline{m}_2'+\frac{a_k(1+c_1m_2+c_1\lambda_k^{\bbC}m_2')}{\lambda_k^{\bbC}(1+c_1m_2)^2}]^2}\nonumber\\
	&\times\left(\underline{m}_2'+\frac{a_k^2c_1m_2'}{(\lambda_k^{\bbC})^2(1+c_1m_2)^4}+\frac{2a_k(1+\underline{m}_2+\lambda_k^{\bbC}\underline{m}_2')}{(\lambda_k^{\bbC})^2(1+c_1m_2)^2}\right),
\end{align}
where $m_2$ and $\underline{m}_2$ are defined in (\ref{m_2}) and (\ref{munderline}), $m_2'$ and $\underline{m}_2'$ stand for their derivatives at $\lambda_k^{\bbC}$, respectively. We shall have establish the proof of theorem \ref{CLTtuS1} if we prove the limiting distribution of $\boldsymbol{\Omega}_0^{\bbC_n}$ and the limits of $\varepsilon_1$, $\varepsilon_2$ and $\varepsilon_3$. In the following, we will put these proofs.

\subsubsection{Limiting distribution of $\boldsymbol{\Omega}_0^{\bbC_n}$}

In this section, we proceed to show the limiting distribution of $\boldsymbol{\Omega}_0^{\bbC_n}$. According to the definition of $\boldsymbol{\Omega}_0^{\bbC_n}$ in (\ref{305}), the proof will be divided into three parts, where are $(\bbX_1\bbA\bbX_1^*\!-\!({\rm tr}\bbA)\bbI_M)/n$,  $(\bbX_1\bbA\boldsymbol\Xi_1^* 
\!+\!\boldsymbol\Xi_1\bbA\bbX_1^*)/n$ and
\begin{align}\label{second2}
	\frac{1}{n}\bbD_{11}^*\left[\bbA_{11}\!+\!\frac{\bbI}{\lambda_{nk}^{\bbC_n}(1+c_{1n}m_{2n}(\lambda_{nk}^{\bbC_n}))}\right]\bbD_{11}.
\end{align}
By Theorem 7.1 in \cite{BY2008}, it is easy to obtain that the $k$-th block of $M\times M$ matrix 
$\frac1{\sqrt{n}}\bbX_1\bbA\bbX_1^*-\frac1{\sqrt{n}}({\rm tr} \bbA(\lambda_{nk}^{\bbC_n}))\bbI_M$ tends to $m_k$-dimensional GOE (GUE) matrix with scale parameter $\underline{m}_2'(\lambda_{k}^{\bbC})$ under real (complex) case.


Having disposed of $\frac1{\sqrt{n}}\bbX_1\bbA\bbX_1^*-\frac1{\sqrt{n}}({\rm tr} \bbA(\lambda_{nk}^{\bbC_n}))\bbI_M$, we can now turn to the proof of (\ref{second2}). To complete the proof of (\ref{second2}), we consider the limiting distribution of 
\begin{align*}
	\bbA_{11}(\lambda_{nk}^{\bbC})+\frac1{\lambda_{nk}^{\bbC}(1+c_{1n}m_{2n}(\lambda_{nk}^{\bbC}))}\bbI_M.
\end{align*}
Rewrite $\boldsymbol\Xi_2$ as $\boldsymbol\Xi_2=(\boldsymbol{0}_{p-M,M},(\tilde{\bbD}_{22})_{p-M,n-M})$
and similar considerations apply to $\bbX_2$,  $\bbX_2=((\bbX_{21})_{p-M,M},(\bbX_{22})_{p-M,n-M})$. Then, we have
\begin{equation}\label{501}
	\bbA=\begin{pmatrix}\frac{1}{n}\bbX_{21}^{\ast}\bbX_{21}\!-\!\lambda_{nk}^{\bbC}\bbI_M, &\!-\!\frac{1}{n}\bbX_{21}^{\ast}(\bbX_{22}\!+\!\tilde{\bbD}_{22})\cr
		\!-\!\frac{1}{n}(\bbX_{22}\!+\!\tilde{\bbD}_{22})^{\ast}\bbX_{21}, &\frac{1}{n}(\bbX_{22}\!+\!\tilde{\bbD}_{22})^{\ast}(\bbX_{22}\!+\!\tilde{\bbD}_{22})\!-\!\lambda_{nk}^{\bbC}\bbI_{n-M}\end{pmatrix}^{-1}.
\end{equation}
Then, the first $M\times M$ major diagonal submatrix is
\begin{align*}
	\bbA_{11}=&\left(\frac{1}{n}\bbX_{21}^{\ast}\bbX_{21}\!-\!\lambda_{nk}^{\bbC}\bbI_M\!-\!\frac{1}{n}\bbX_{21}^{\ast}(\bbX_{22}\!+\!\tilde{\bbD}_{22})\bbA_{22}\frac{1}{n}(\bbX_{22}\!+\!\tilde{\bbD}_{22})^{\ast}\bbX_{21}\right)^{-1}\\
	=&\left(\!-\!\lambda_{nk}^{\bbC}\bbI_M\!+\!\frac{1}{n}\bbX_{21}^{\ast}\left[\bbI_{p\!-\!M}\!-\!\frac{1}{n}(\tilde{\bbD}_{22}\!+\!\bbX_{22})\bbA_{22}(\tilde{\bbD}_{22}\!+\!\bbX_{22})^{\ast}\right]\bbX_{21}\right)^{-1}\\
	=&\left(\!-\!\lambda_{nk}^{\bbC}\bbI_M\!-\!\frac{\lambda_{nk}^{\bbC}}{n}\bbX_{21}^{\ast}\tilde{\bbA}_{22}\bbX_{21}\right)^{-1}\\
	=&\left(\!-\!\lambda_{nk}^{\bbC}\bbI_M\!-\!\frac{\lambda_{nk}^{\bbC}}{n}(\tr\tilde{\bbA}_{22})\bbI_M\!-\!\boldsymbol{\Omega}_1^{\bbC_n}(\lambda_{nk}^{\bbC},\bbX_{21})\right)^{-1},
\end{align*}
where 
\begin{align*}
	&\bbA_{22}=[\frac{1}{n}(\bbX_{22}+\tilde{\bbD}_{22})^{\ast}(\bbX_{22}+\tilde{\bbD}_{22})-\lambda_{nk}^{\bbC}\bbI_{n-M}]^{-1},\\
	&\tilde{\bbA}_{22}=[\frac{1}{n}(\bbX_{22}+\tilde{\bbD}_{22})(\bbX_{22}+\tilde{\bbD}_{22})^{\ast}-\lambda_{nk}^{\bbC}\bbI_{p-M}]^{-1},\\
	&\boldsymbol{\Omega}_1^{\bbC_n}(\lambda_{nk}^{\bbC},\bbX_{21})=\frac{\lambda_{nk}^{\bbC}}{n}\left[\bbX_{21}^{\ast}\tilde{\bbA}_{22}\bbX_{21}-(\tr\tilde{\bbA}_{22})\bbI_M\right].
\end{align*}

We emphasize that both $\bbA_{22}$ and $\bbA$ are noncentral sample covariance matrices with the same limiting noncentral parameter matrix. So the Stieltjes transform of LSD of $\bbA_{22}$ is $\underline{m}_2(\cdot)$. Similar arguments apply to $\tilde{\bbA}_{22}$, we conclude that  $\frac{1}{p-M}\tr\tilde{\bbA}_{22}$ tends to $m_2(\lambda_k^{\bbC})$ with probability one. By the CLT of the quadratic form, we have
\begin{equation*}
	[\boldsymbol{\Omega}_1^{\bbC_n}(\lambda_{nk}^{\bbC},\bbX_{21})]_{ij}=O_p\left(\frac{\lambda_{nk}^{\bbC}}{\sqrt{n}}\right),
\end{equation*}
and the $k$-th block of $\frac{1}{\sqrt{p-M}}(\bbX_{21}^*\tilde{\bbA}_{22}\bbX_{21}-(\tr\tilde{\bbA}_{22})\bbI_{M})$ tends to $m_k$-dimensional GOE (GUE) matrix with scale parameter $m_2'(\lambda_{k}^{\bbC})$ under real (complex) case. 
Moreover,
\begin{eqnarray}\label{502}
	&&\bbA_{11}-\frac{-1}{\lambda_{nk}^{\bbC}+\lambda_{nk}^{\bbC}\frac{p-M}{n}m_{2n}(\lambda_{nk}^{\bbC})}\bbI_M\nonumber\\
	=&&\frac{-1}{\lambda_{nk}^{\bbC}+\lambda_{nk}^{\bbC}\frac{p-M}{n}m_{2n}(\lambda_{nk}^{\bbC})}\boldsymbol{\Omega}_1^{\bbC_n}(\lambda_{nk}^{\bbC},\bbX_{21})\bbA_{11}\nonumber\\
	=&&\frac{\boldsymbol{\Omega}_1^{\bbC_n}(\lambda_{nk}^{\bbC},\bbX_{21})}{[\lambda_{nk}^{\bbC}+\lambda_{nk}^{\bbC}\frac{p-M}{n}m_{2n}(\lambda_{nk}^{\bbC})]^2}+\frac{[\boldsymbol{\Omega}_1^{\bbC_n}(\lambda_{nk}^{\bbC},\bbX_{21})]^2\bbA_{11}}{[\lambda_{nk}^{\bbC}+\lambda_{nk}^{\bbC}\frac{p-M}{n}m_{2n}(\lambda_{nk}^{\bbC})]^2}
\end{eqnarray}
where 
\begin{align*}
	\frac{[\boldsymbol{\Omega}_1^{\bbC_n}(\lambda_{nk}^{\bbC},\bbX_{21})]^2\bbA_{11}}{[\lambda_{nk}^{\bbC}+\lambda_{nk}^{\bbC}\frac{p-M}{n}m_{2n}(\lambda_{nk}^{\bbC})]^2}=O_p\left(\frac{1}{n}\right)\|\bbA_{11}\|\boldsymbol{1}\boldsymbol{1}'.
\end{align*}
By (\ref{502}) and classical CLT, it is easy to see that the corresponding block of $\frac{1}{\sqrt{n}}\bbD_{11}[\bbA_{11}+\frac1{l_i^{\bbC_n}(1+c_{3n}m_{2n}(l_i^{\bbC_n}))}\bbI_M]\bbD_{11}$ tends to the GOE (GUE) matrix under real (complex) case with scale parameter 
\begin{align*}
	\frac{a_k^2c_1m_2'(\lambda_k^{\bbC})}{(\lambda_k^{\bbC})^2(1+c_1m_2(\lambda_k^{\bbC}))^4}.
\end{align*}

We shall have established the limiting distribution of 
$\boldsymbol{\Omega}_0^{\bbC_n}$ if we give the limiting distribution of
\begin{align}\label{third3}
	\frac{1}{\sqrt{n}}\bbX_1\bbA\boldsymbol\Xi_1^*+\frac{1}{\sqrt{n}}\boldsymbol\Xi_1\bbA\bbX_1^*.
\end{align}
It is easily seen that the elements $(s,t)$ ($s\leq M, t\leq M$) of (\ref{third3})  is $1/\sqrt{n}d_t\bbx_s\bba_t+1/\sqrt{n}d_s\bba_s^*\bbx_t^*$, where $\bbx_s$ is the $s$-th row of $\bbX_1$ and $\bba_t$ is the $t$-th column of $\bbA$. Since $\bbX_1$ is independent of $\bbA$,  the limiting distribution of $1/\sqrt{n}d_t\bbx_s\bba_t+1/\sqrt{n}d_s\bba_s^*\bbx_t^*$ is Gaussian under given $\bbA$. The mean of the Gaussian is zero and its variance is present under two different situations. Under the real samples situation, the variance is equal to 
\begin{equation}\label{504}
	\sigma_n^2(s,t)=\begin{cases}4a_s\E\bba_s^{T}\bba_s,& \mbox{ if } s=t\cr 
		a_s\E\bba_s^T\bba_s+a_t\E\bba_t^T\bba_t,&\mbox{ if } s\ne t\cr\end{cases}
\end{equation}
and	under the complex samples situation, the variance is equal to
{\begin{equation}\label{505}
		\sigma_n^2(s,t)=a_s\E\bba_s^*\bba_s+a_t\E\bba_t^*\bba_t.
\end{equation}}
What is left is to compute the limits of (\ref{504}) and (\ref{505}). By the inverse blockwise matrix formula, we know that the vector consisted of the first $M$ components of $\bba_t$ is equal to the $t$-th column ${\mathbf A}_{11t}$ of $\bbA_{11}$ and the vector consisted of next $n-M$ components is equal to the $j$-th column of $\bbA_{22}\frac{1}{n}(\bbX_{22}^*+\tilde{\bbD}_{22}^*)\bbX_{21}{\mathbf A}_{11t}$. By (\ref{501}), we have
\begin{align*}
	\E\bba_t^*\bba_t=&\E\bbA_{11t}^*\bbA_{11t}+\E\bbA_{11t}^*\bbX_{21}^*\frac{1}{n}(\bbX_{22}+\tilde{\bbD}_{22})\bbA_{22}^2\frac{1}{n}(\bbX_{22}+\tilde{\bbD}_{22})^*\bbX_{21}\bbA_{11t}\\
	=&\E\bbA_{11t}^*\left[\bbI_M+\frac{1}{n}\bbX_{21}^*\frac{1}{n}(\bbX_{22}+\tilde{\bbD}_{22})\bbA_{22}^2(\bbX_{22}^*+\tilde{\bbD}_{22}^*)\bbX_{21}\right]\bbA_{11t}\\
	=\E\bbA_{11t}^*&\left[\bbI_M+\frac{1}{n}\tr\left[\frac{1}{n}(\bbX_{22}+\tilde{\bbD}_{22})\bbA_{22}^2(\bbX_{22}+\bbD_{22})^*\right]{\bbI}_M+O_p\left(\frac{1}{\sqrt{n}}\right)\right]\bbA_{11t},
\end{align*}
where
\begin{align*}
	&\left(1+\frac{1}{n}\tr\frac{1}{n}(\bbX_{22}+\tilde{\bbD}_{22})\bbA_{22}^2(\bbX_{22}+\tilde{\bbD}_{22})^*\right)\\
	=&\left(1+\frac{1}{n}\tr\bbA_{22}^2\frac{1}{n}(\bbX_{22}+\tilde{\bbD}_{22})^*(\bbX_{22}+\tilde{\bbD}_{22})\right)\\
	=&\left(1+\frac{1}{n}\tr\bbA_{22}+\frac{\lambda_{nk}^{\bbC_n}}{n}\tr\bbA_{22}^2\right)=1+\underline m_{2n}(\lambda_{nk}^{\bbC_n})+\lambda_{nk}^{\bbC_n}\underline m_{2n}'(\lambda_{nk}^{\bbC_n}).
\end{align*}
Then applying (\ref{502}), we have
\begin{align}\label{eqaa1}
	\E\bba_t^*\bba_t=&\frac{\left(1+\underline m_{2n}(\lambda_{nk}^{\bbC_n})+\lambda_{nk}^{\bbC_n}\underline m_{2n}'(\lambda_{nk}^{\bbC_n})\right)}{(\lambda_{nk}^{\bbC_n})^2(1+c_{1n}m_{2n}(\lambda_{nk}^{\bbC_n}))^2}\left(1+o\left(\frac{1}{n}\right)\right)\nonumber\\
	\to&\frac{1}{(\lambda_k^{\bbC})^2(1+c_1m_2(\lambda_k^{\bbC}))^2}\left(1+\underline m_2(\lambda_k^{\bbC})+\lambda_k^{\bbC}\underline m_2'(\lambda_k^{\bbC})\right),
\end{align}
In the same manner we can see that 
\begin{equation}\label{eqaij}
	\E\bba_s^*\bba_t\to 0, \mbox{ for } s\ne t,
\end{equation}
which implies the variables in asymmetric positions are asymptotic independent.
From the above analysis, we obtain the corresponding block of $\frac{1}{\sqrt{n}}\bbX_1\bbA\boldsymbol\Xi_1^*+\frac{1}{\sqrt{n}}\boldsymbol\Xi_1\bbA\bbX_1^*$ tends to a $m_k\times m_k$ GOE (GUE) matrix under real (complex) case with scale parameter 
\begin{align*}
	\frac{2a_k(1+\underline{m}_2(\lambda_k^{\bbC})+\lambda_k^{\bbC}\underline{m}_2'(\lambda_k^{\bbC}))}{(\lambda_k^{\bbC})^2(1+c_1m_2(\lambda_k^{\bbC}))^2}.
\end{align*}

The third moment of normal population is zero so that the limiting distribution of $\frac1n\bbX_1\bbA\bbX_1^*\!-\!\frac1n({\rm tr}\bbA)\bbI_M$,  $\frac1n(\bbX_1\bbA\boldsymbol\Xi_1^* 
\!+\!\boldsymbol\Xi_1\bbA\bbX_1^*)$ and (\ref{second2}) are independent among them.
We conclude that the limiting distribution of the $k$-th block of $\boldsymbol{\Omega}_0^{\bbC_n}$ equals the $m_k\times m_k$ GOE (GUE) matrix under real (complex) matrix with the scale parameter  
\begin{align}
	\underline{m}_2'(\lambda_k^{\bbC})+\frac{a_k^2c_1m_2'(\lambda_k^{\bbC})}{(\lambda_k^{\bbC})^2(1+c_1m_2(\lambda_k^{\bbC}))^4}+\frac{2a_k(1+\underline{m}_2(\lambda_k^{\bbC})+\lambda_k^{\bbC}\underline{m}_2'(\lambda_k^{\bbC}))}{(\lambda_k^{\bbC})^2(1+c_1m_2(\lambda_k^{\bbC}))^2}.
\end{align}

\subsubsection{Limits of $\varepsilon_1$, $\varepsilon_2$ and $\varepsilon_3$}
In this section, we proceed to show the limits of $\varepsilon_1$, $\varepsilon_2$ and $\varepsilon_3$ defined in (\ref{e1})-(\ref{e3}). To deal with  $\varepsilon_1$,  we note that
\begin{align*}
	\varepsilon_1&\!=\!\frac{\gamma_{k}^{\bbC_n}}{\sqrt{n}}\frac{\lambda_{nk}^{\bbC}}{n}\tr\bbA_n(\lambda_{nk}^{\bbC})\bbA_n(l_i^{\bbC_n})+\frac{1}{n}\tr\bbA_n(\lambda_{nk}^{\bbC})-\underline{m}_{2n}^0(\lambda_{nk}^{\bbC})\\
	&=\frac{\gamma_{k}^{\bbC_n}}{\sqrt{n}}\frac{\lambda_{nk}^{\bbC}}{n}\tr\bbA_n^2(\lambda_{nk}^{\bbC})\!+\!\frac{\gamma_{k}^{\bbC_n}}{\sqrt{n}}\frac{\lambda_{nk}^{\bbC}}{n}\tr\bbA_n(\lambda_{nk}^{\bbC})[\bbA_n(l_i^{\bbC_n})\!-\!\bbA_n(\lambda_{nk}^{\bbC})]+O_p\left(\frac{1}{n}\right),
\end{align*}
where the last equality is the consequence of Theorem 1 in \cite{BannaN20CLT}.
As $\|\bbA_n(l_i^{\bbC_n})\|$ is bounded a.s. and Theorem \ref{limittuS1}, we have 
\begin{align*}
	&\frac{1}{n}\tr\bbA_n^2(\lambda_{nk}^{\bbC})\overset{a.s.}{\longrightarrow}\underline{m}_2'(\lambda_{nk}^{\bbC}),\\
	&\frac{\gamma_{k}^{\bbC_n}}{\sqrt{n}}\frac{\lambda_{nk}^{\bbC}}{n}\tr\bbA_n(\lambda_{nk}^{\bbC})\left[\bbA_n(l_i^{\bbC_n})\!-\!\bbA_n(\lambda_{nk}^{\bbC})\right]=o_p\left(\frac{\gamma_k^{\bbC_n}}{\sqrt{n}}\lambda_{nk}^{\bbC}\right).
\end{align*}

The task is now to consider the limit of $\varepsilon_2$,
\begin{align*}
	\varepsilon_2=&\frac{-1}{l_i^{\bbC_n}+l_i^{\bbC_n}\frac{p-M}{n}m_{2n}(l_i^{\bbC_n})}-\frac{-1}{\lambda_{nk}^{\bbC}+\lambda_{nk}^{\bbC}c_{1n}m_{2n}^0(\lambda_{nk}^{\bbC})}\\
	=&\frac{l_i^{\bbC_n}(1+\frac{p-M}{n}m_{2n}(l_i^{\bbC_n}))-\lambda_{nk}^{\bbC}(1+c_{1n}m_{2n}^0(\lambda_{nk}^{\bbC}))}{[\lambda_{nk}^{\bbC}+\lambda_{nk}^{\bbC}c_{1n}m_{2n}^0(\lambda_{nk}^{\bbC})]^2}\\
	+&\frac{(l_i^{\bbC_n}-\lambda_{nk}^{\bbC})^2(1+c_{1n}m_{2n}(l_i^{\bbC_n})+\lambda_{nk}^{\bbC}c_{1n}(m_{2n}^0)'(\lambda_{nk}^{\bbC}))^2+o_p(\frac{1}{n})}{[\lambda_{nk}^{\bbC}+\lambda_{nk}^{\bbC}c_{1n}m_{2n}^0(\lambda_{nk}^{\bbC})]^2(l_i^{\bbC_n}+l_i^{\bbC_n}\frac{p-M}{n}m_{2n}^0(l_i^{\bbC_n}))}\\
	=&\frac{\gamma_{k}^{\bbC_n}}{\sqrt{n}}\left[\frac{1+c_1m_2(\lambda_{k}^{\bbC})+c_1\lambda_{k}^{\bbC}m_2'(\lambda_{k}^{\bbC})}{\lambda_{k}^{\bbC}[1+c_1m_2(\lambda_{k}^{\bbC})]^2}+o_p(1)\right],
\end{align*}
the last equality being a consequence of 
\begin{align*}
	\frac{(l_i^{\bbC_n}-\lambda_{nk}^{\bbC})^2(1+c_{1n}m_{2n}(l_i^{\bbC_n})+\lambda_{nk}^{\bbC}c_{1n}(m_{2n}^0)'(\lambda_{nk}^{\bbC}))^2+o_p(\frac{1}{n})}{[\lambda_{nk}^{\bbC}+\lambda_{nk}^{\bbC}c_{1n}m_{2n}^0(\lambda_{nk}^{\bbC})]^2(l_i^{\bbC_n}+l_i^{\bbC_n}\frac{p-M}{n}m_{2n}(l_i^{\bbC_n}))}=o_p\left(\frac{\gamma_i^{\bbC_n}}{\sqrt{n}}\right).
\end{align*}

It remains to show the limit of $\varepsilon_3$. Here, we focus on  $\frac{1}{n}\bbX_1\bbA\bbX_1^*-\frac{1}{n}(\tr\bbA)\bbI_M$ and other terms can be handle in a similar way. According to $a_k$, we have
\begin{align*}
	&\frac{1}{n}\bbX_1\bbA(l_i^{\bbC_n})\bbX_1^*-\frac{1}{n}(\tr\bbA(l_i^{\bbC_n}))\bbI_M-(\frac{1}{n}\bbX_1\bbA(\lambda_{nk}^{\bbC})\bbX_1-\frac{1}{n}(\tr\bbA(\lambda_{nk}^{\bbC}))\bbI_M)\\
	=&\left[\frac{1}{n}\bbX_1\bbA(\lambda_{nk}^{\bbC})\bbA(l_i^{\bbC_n})\bbX_1^*-\frac{1}{n}\tr\bbA(\lambda_{nk}^{\bbC})\bbA(l_i^{\bbC_n})\bbI_M\right](l_i^{\bbC_n}-\lambda_{nk}^{\bbC})\\
	=&O_p\left(\frac{\lambda_k^{\bbC}}{\sqrt{n}}\right)\frac{\gamma_k^{\bbC_n}}{\sqrt{n}}\boldsymbol{1}\boldsymbol{1}'=o_p\left(\frac{\gamma_k^{\bbC_n}}{\sqrt{n}}\right)\boldsymbol{1}\boldsymbol{1}'.
\end{align*}
Note that the assumption the rate of $a_k$ diverging to infinity in  Assumption A is used here. To be specific,  We find that if the rate of $a_k$ diverging to infinity is more than $\sqrt{n}$, ${\lambda_k^{\bbC}}/{\sqrt{n}}$ will tend infinity.

Combining the limiting distribution of $\boldsymbol{\Omega}_0^{\bbC_n}$, we obtain the limiting distribution of the random vector $\{\gamma_k^{\bbC_n}\}$ and its key scale parameter is defined in (\ref{theta1}).

\subsection{Proof of Theorem \ref{limitsF}}
\label{Phase T3}
At first, we consider the first order limit of the spiked eigenvalues of $\mathbf{F}_p$ under given the matrix sequence $\{\bbC_n\}$.  Recall 
\begin{eqnarray*}
	\bbC_n&=&\frac{1}{n}(\boldsymbol \Xi+\bbX)(\boldsymbol \Xi+\bbX)^{\ast}=\bbU\left(\begin{array}{cc}{\boldsymbol\Sigma}_1&0\\0&{\boldsymbol\Sigma}_2\end{array}\right)\bbU^{\ast},\\
	\bbS_N&=&\frac{1}{N}\bbY_N\bbY_N^{\ast}=\frac{1}{N}\left(\begin{array}{c}\bbY_1\\\bbY_2\end{array}\right)\left(\begin{array}{cc}\bbY_1^{\ast}&\bbY_2^{\ast}\end{array}\right)=\frac{1}{N}\left(\begin{array}{cc}\bbY_1\bbY_1^{\ast} &\bbY_1\bbY_2^{\ast}\\ \bbY_2\bbY_1^{\ast} & \bbY_2\bbY_2^{\ast} \end{array}\right),
\end{eqnarray*}
where ${\boldsymbol\Sigma}_1$ is an $M\times M$ diagonal matrix and $\bbY_1$ denotes the first $M$ rows of $\bbY$.
The matrix $\bbC_n$ can be seen as a general non-negative definite matrix with the eigenvalues formed in descending order,
\begin{align}\label{eigtuS1inF}
	l_1^{\bbC_n}\geq l_2^{\bbC_n} \geq \cdots \geq l_p^{\bbC_n}.
\end{align}
We consider the eigen-equation
\begin{align*}
	&\left|\bbC_n\bbS_N^{-1}-\lambda\bbI\right|=0\Longleftrightarrow\left|	\bbC_n-\lambda\bbS_N\right|=0\\
	\Longleftrightarrow&\left|\bbU\left(\begin{array}{cc}\boldsymbol\Sigma_1&0\\0&\boldsymbol\Sigma_2\end{array}\right)\bbU^{\ast}-\frac{\lambda}{N}\bbY_N\bbY_N^{\ast}\right|=0\\
	\Longleftrightarrow&\left|\left(\begin{array}{cc}\boldsymbol\Sigma_1&0\\0&\boldsymbol\Sigma_2\end{array}\right)-\frac{\lambda}{N}\bbU^{\ast}\bbY_N\bbY_N^{\ast}\bbU\right|=0.
\end{align*}

We denote $\bbU^{\ast}\bbY_N\bbY_N^{\ast}\bbU/N$ as $\tilde{\bbY}_N\tilde{\bbY}_N^{\ast}/N$. The entries of $\bbY_N$ are the standard normal so that both  $\tilde{\bbY}_N\tilde{\bbY}_N^{\ast}/N$ and ${\bbY}_N{\bbY}_N^{\ast}/N$ have the same distribution. If there is no confusion, we will still write the notation $\bbY_N$. Then the eigen-equation becomes 
\begin{align*}
	&\left|\left(\begin{array}{cc}\boldsymbol\Sigma_1&0\\0&\boldsymbol\Sigma_2\end{array}\right)-\frac{\lambda}{N}\bbY_N\bbY_N^{\ast}\right|=0\\
	\Longleftrightarrow&\left|\left(\begin{array}{cc}\boldsymbol\Sigma_1&0\\0&\boldsymbol\Sigma_2\end{array}\right)-\lambda\left(\begin{array}{cc}\frac{1}{N}\bbY_1\bbY_1^{\ast}&\frac{1}{N}\bbY_1\bbY_2^{\ast}\\\frac{1}{N}\bbY_2\bbY_1^{\ast}&\frac{1}{N}\bbY_2\bbY_2^{\ast}\end{array}\right)\right|=0.
\end{align*}
According to the sample spiked eigenvalues $l_i$, $i\in\mathcal{J}_k$, $k=1,\cdots,K$ of $\boldsymbol F_p=\bbC_n\bbS_N^{-1}$,  we have $|\boldsymbol\Sigma_2-l_i\frac{1}{N}\bbY_2\bbY_2^{\ast}|\neq 0$ almost surely, then 
\begin{eqnarray*}
	&&\left|\boldsymbol\Sigma_1-l_i\frac{1}{N}\bbY_1\bbY_1^{\ast}-l_i^2\frac{1}{N}\bbY_1\bbY_2^{\ast}\left(\boldsymbol\Sigma_2-l_i\frac{1}{N}\bbY_2\bbY_2^{\ast}\right)^{-1}\frac{1}{N}\bbY_2\bbY_1^{\ast}\right|=0\\
	&&\left|\boldsymbol\Sigma_1-l_i\frac{1}{N}\bbY_1\left[\bbI_N+l_i\frac{1}{N}\bbY_2^{\ast}(\boldsymbol\Sigma_2-l_i\frac{1}{N}\bbY_2\bbY_2^{\ast})^{-1}\bbY_2\right]\bbY_1^{\ast}\right|=0\\
	&&\left|\boldsymbol\Sigma_1-l_i\frac{1}{N}\tr\left[\bbI_N+l_i\frac{1}{N}\bbY_2^{\ast}(\boldsymbol\Sigma_2-l_i\frac{1}{N}\bbY_2\bbY_2^{\ast})^{-1}\bbY_2\right]\bbI_M+\boldsymbol{\Omega}_q^{\boldsymbol F}(l_i)\right|=0\\
	&&\left|\boldsymbol\Sigma_1\!-\!l_i\left(1\!+\!l_i\frac{p-M}{N}\frac{1}{p-M}\tr\left[(\boldsymbol\Sigma_2\!-\!l_i\frac{1}{N}\bbY_2\bbY_2^{\ast})^{-1}\frac{1}{N}\bbY_2\bbY_2^{\ast}\right]\right)\bbI_M\!+\!\boldsymbol{\Omega}_q^{\boldsymbol F}(l_i)\right|=0\\
	&&\left|\boldsymbol\Sigma_1\!-\!l_i\left(1\!+\!l_i\frac{p-M}{N}\frac{1}{p-M}\tr\left(\boldsymbol\Sigma_2(\frac{1}{N}\bbY_2\bbY_2^{\ast})^{-1}\!-\!l_i\bbI\right)^{-1}\right)\bbI+\boldsymbol{\Omega}_q^{\boldsymbol F}(l_i)\right|=0,
\end{eqnarray*}
where 
\begin{eqnarray}\label{307}
	\boldsymbol{\Omega}_N^{\boldsymbol F}(l_i)=&&\frac{l_i}{N}\bbY_1\left[\bbI_N+l_i\frac{1}{N}\bbY_2^{\ast}\left(\boldsymbol\Sigma_2-l_i\frac{1}{N}\bbY_2\bbY_2^{\ast}\right)^{-1}\bbY_2\right]\bbY_1^{\ast}\\
	&&-\frac{l_i}{N}\tr\left[\bbI_N+l_i\frac{1}{N}\bbY_2^{\ast}\left(\boldsymbol \Sigma_2-l_i\frac{1}{N}\bbY_2\bbY_2^{\ast}\right)^{-1}\bbY_2\right]{\bbI_M}.
\end{eqnarray}
Similar to Subsection \ref{523}, by the assumption on $a_k$, the strong law of large numbers and the Stietjes transformation equation of $M$-$P$ law, we have $\boldsymbol{\Omega}_N^{\boldsymbol F}(l_i)\overset{a.s.}{\longrightarrow}\boldsymbol{0}_{M\times M}$ as $N\to\infty$, and 
\begin{eqnarray}
	&\frac{1}{p-M}\tr\left(\boldsymbol \Sigma_2(\frac{1}{N}\bbY_2\bbY_2^{\ast})^{-1}-z\bbI\right)^{-1}\overset{a.s.}{\longrightarrow}m_{3}(z),\quad {\forall z\in\mathbb{C}^{+}\cup S_F^c} \label{3042}
\end{eqnarray}
where $m_{3}(\cdot)$ is the Stieltjes transform  of $F$. 
Taking limit on the eigen-equation, there exists a diagonal block being $\boldsymbol{0}_{m_i\times m_i}$, then we have 
$\lambda_k^{\bbC}=\lambda_k(1+c_2\lambda_km_{3}(\lambda_k))$, where $\lambda_k^{\bbC}$ and $\lambda_k$ are the almost surely limit of $l_i^{\bbC_n}$ and $l_i, i\in \mathcal{J}_k$ respectively.
According to (\ref{m_3}), we have 
\begin{eqnarray*}
	zm_{3}(z)=&&(z+c_2z^2m_{3}(z))m_2(z+c_2z^2m_3(z))\\
	\Longleftrightarrow \lambda_km_{3}(\lambda_k)=&&\lambda_k(1+c_2\lambda_km_{3}(\lambda_k))m_2(\lambda_k(1+c_2\lambda_km_{3}(\lambda_k)))=\lambda_k^{\bbC}m_2(\lambda_k^{\bbC}),
\end{eqnarray*} 	
i.e., 
\begin{eqnarray}
	\lambda_k^{\bbC}&=&\lambda_k(1+c_2\lambda_k^{\bbC}m_2(\lambda_k^{\bbC}))\label{relation3}\\
	\psi_{\boldsymbol F}(\lambda_k^{\bbC})&\overset{\bigtriangleup}{=}&\lambda_k=\frac{\lambda_k^{\bbC}}{1+c_2\lambda_k^{\bbC}m_2(\lambda_k^{\bbC})}=\frac{\lambda_k^{\bbC}}{1+c_2\lambda_km_3(\lambda_k)}.\label{relation302}
\end{eqnarray}

From what has already been proved, we conclude that the first order limit of $l_i$ is independent of $\{\bbC_n\}$, only related to the limit of their spiked eigenvalues. An easy computation shows the relationship between $\lambda_k$ and $a_k$:
\begin{align*}
	\lambda_k=\frac{\lambda_k^{\bbC}}{1+c_2\lambda_k^{\bbC}m_2(\lambda_k^{\bbC})},\quad \lambda^{\bbC}_k=a_k\left(1-c_1m_1(a_k)\right)^2+(1-c_1)\left(1-c_1m_1(a_k)\right).
\end{align*}

\subsection{Proof of Theorem \ref{CLTnonF}}\label{CLT3}
We are now in a position to study the asymptotic distribution of the random vector 
\begin{align}\label{gama}
	\gamma_{Nk}=(\sqrt{N}(l_i-\lambda_{Nk})/\lambda_{Nk}, i\in\mathcal{J}_k)
\end{align}
{where $\lambda_{Nk}=\psi_{\bbF}(\lambda_{nk}^{\bbC_n})$ with the parameters in function $\psi_{\bbF}$ replaced by the corresponding empirical parameters, and} which will be divided into two steps. At first, we give the conditional limiting distribution of $\gamma_{Nk}|\bbC_n$, then we will find  the limiting distribution are independent with the {choice of} conditioning $\bbC_n$.  Secondly, combining the above theorems and the following subsection, we can complete the proof of Theorem \ref{CLTnonF}.
\subsubsection{The conditional limiting distribution of $\gamma_{Nk}|\bbC_n$}
In this section, we consider the CLT of the random vector $\gamma_{Nk}|\bbC_n=\{\sqrt{N}(l_i-\psi_{\boldsymbol{F}}(l_i^{\bbC_n}))/\psi_{\boldsymbol{F}}(l_i^{\bbC_n}), i\in\mathcal{J}_k\}$.
Recall the eigen-equation 
\begin{eqnarray*}
	&\left|\boldsymbol\Sigma_1\!-\!l_i\left(1\!+\!l_i\frac{p\!-\!M}{N}\frac{1}{p\!-\!M}\tr\left(\boldsymbol\Sigma_2\left(\frac{1}{N}\bbY_2\bbY_2^{\ast}\right)^{-1}\!-\!l_i\bbI_{p-M}\right)^{-1}\right)\bbI_M\!+\!\boldsymbol{\Omega}_N^{\boldsymbol F}(l_i)\right|\!=\!0\\
	&\left|\boldsymbol\Sigma_1\!-\!\psi_{\boldsymbol{F}}(l_i^{\bbC_n})\left(1\!+\!\frac{p\!-\!M}{N}\psi_{\boldsymbol{F}}(l_i^{\bbC_n})m_{3N}\psi_{\boldsymbol{F}}(l_i^{\bbC_n})\right)\bbI_M\!+\!\boldsymbol{\Omega}_N^{\boldsymbol F}(\psi_{\boldsymbol{F}}(l_i^{\bbC_n}))\!+\!\varepsilon_1\right|=0
\end{eqnarray*}
where 
\begin{eqnarray*}
	&&m_{3N}(\psi_{\boldsymbol{F}}(l_i^{\bbC_n}))=\frac{1}{p\!-\!M}\tr\left(\boldsymbol\Sigma_2\left(\frac{1}{N}\bbY_2\bbY_2^{\ast}\right)^{\!-\!1}\!-\!\psi_{\boldsymbol{F}}(l_i^{\bbC_n})\bbI_{p-M}\right)^{-1}\\
	&&\varepsilon_1=\psi_{\boldsymbol{F}}(l_i^{\bbC_n})\left(1+\psi_{\boldsymbol{F}}(l_i^{\bbC_n})\frac{p-M}{N}m_{3N}(\psi_{\boldsymbol{F}}(l_i^{\bbC_n}))\right)\bbI_M\\
	&&\quad\quad-l_i\left(1+l_i\frac{p-M}{N}m_{3N}(l_i)\right)\bbI_M+\boldsymbol{\Omega}_N^{\boldsymbol F}(l_i)-\boldsymbol{\Omega}_N^{\boldsymbol F}(\psi_{\boldsymbol{F}}(l_i^{\bbC_n})).
\end{eqnarray*}
By simple calculation, we set $c_{2N}=(p-M)/N$ and have  
\begin{align*}
	\varepsilon_1=&[-(l_i-\psi_{\boldsymbol{F}}(l_i^{\bbC_n}))-(l_i^2-(\psi_{\boldsymbol{F}}(\l_i^{\bbC_n}))^2)c_{2N}m_{3N}(l_i)\\
	&-(\psi_{\boldsymbol{F}}(l_i^{\bbC_n}))^2c_{2N}m_{3	N}'(\psi_{\boldsymbol{F}}(l_i^{\bbC_n}))(l_i-\psi_{\boldsymbol{F}}(l_i^{\bbC_n}))(1+o(1))]\mathbf I_M\\
	&+\boldsymbol{\Omega}_N^{\boldsymbol F}(l_i)-\boldsymbol{\Omega}_N^{\boldsymbol F}(\psi_{\boldsymbol{F}}(l_i^{\bbC_n}))
\end{align*}
i.e.
\begin{eqnarray}\label{404}
	\varepsilon_1\!=\!-\!\frac{\gamma_{Nk}|\bbC_n}{\sqrt{N}}\lambda_k\left[1\!+\!2\lambda_kc_{2}m_3(\lambda_k)\!+\!\lambda_k^2c_{2}m_3'(\lambda_k)\right](1\!+\!o_p(1))\bbI\!+\!o_p(\frac{l_i}{\sqrt{N}})\boldsymbol{1}\boldsymbol{1}'.
\end{eqnarray}
By (\ref{relation302}), we have 
\begin{equation*}
		l_i^{\bbC_n}=\psi_F(l_i^{\bbC_n})(1+c_2\psi_{F}(l_i^{\bbC_n})m_3(l_i^{\bbC_n})).
\end{equation*}
We recall 
\begin{align*}
	\left|\boldsymbol\Sigma_1\!-\!\psi_{\boldsymbol{F}}(l_i^{\bbC_n})\left(1\!+\!\psi_{\boldsymbol{F}}(l_i^{\bbC_n})c_{2N}m_{3N}(\psi_{\boldsymbol{F}}(l_i^{\bbC_n}))\right)\bbI_M\!+\!\boldsymbol{\Omega}_N^{\bbF}(\psi_{\boldsymbol{F}}(l_i^{\bbC_n}))+\varepsilon_1\bbI_M\right|=0
\end{align*}
becomes
\begin{align}\label{deter1}
	0=\begin{vmatrix}
		\lambda_1^{\bbC}\!-\!l_i^{\bbC_n}\!+\!O(\frac{\psi_{\boldsymbol{F}}(l_i^{\bbC_n})}{\sqrt{N}})&\cdots&O(\frac{\psi_{\boldsymbol{F}}(l_i^{\bbC_n})}{\sqrt{N}})&\cdots&O(\frac{\psi_{\boldsymbol{F}}(l_i^{\bbC_n})}{\sqrt{N}})\cr
		O(\frac{\psi_{\boldsymbol{F}}(l_i^{\bbC_n})}{\sqrt{N}})&\cdots&\cdots&\cdots&O(\frac{\psi_{\boldsymbol{F}}(l_i^{\bbC_n})}{\sqrt{N}})\cr
		\cdots&\cdots&[\boldsymbol{\Omega}_N^{\boldsymbol F}]_{kk}\!+\!\varepsilon_1\bbI_{m_k}&\cdots&\cdots\cr
		O(\frac{\psi_{\boldsymbol{F}}(l_i^{\bbC_n})}{\sqrt{N}})&\cdots&\cdots&\cdots&O(\frac{\psi_{\boldsymbol{F}}(l_i^{\bbC_n})}{\sqrt{N}})\cr
		O(\frac{\psi_{\boldsymbol{F}}(l_i^{\bbC_n})}{\sqrt{N}})&\cdots&\cdots&\cdots&\lambda_M^{\bbC}\!-\!l_i^{\bbC_n}\!+\!O(\frac{\psi_{\boldsymbol{F}}(l_i^{\bbC_n})}{\sqrt{N}})\cr\end{vmatrix}
\end{align}
where $[\boldsymbol{\Omega}_N^{\boldsymbol{F}}]_{kk}$ is $k$-th diagonal block of $\boldsymbol{\Omega}_N^{\boldsymbol{F}}(\psi_{\boldsymbol{F}}(l_i^{\bbC_n}))$.

By Skorokhod strong representation theorem (for more detais, see \cite{Sk1956} or \cite{HB2014}), on an appropriate 
probability space, one may redefine the random variables such that $\boldsymbol{\Omega}_N^{\boldsymbol{F}}$ tends to the Gaussian variables with probability one.  Multiplying $(\psi_{\boldsymbol{F}}(l_i^{\bbC_n}))^{-1/2}N^{1/4}$ to the $k$-th block row and column of the determinant in (\ref{deter1}),  $p\to\infty$ and $N\to\infty$. {It is easily seen that  all non-diagonal elements tend to zero and and all the diagonal entries except the $k$-th are bounded away from zero as $p\to\infty$ or $N\to\infty$.} Therefore, 
\begin{equation*}
	[\sqrt{N}\boldsymbol{\Omega}_q^{\boldsymbol F}]_{kk}-(\gamma_{Nk}|\bbC_n)\lambda_k\vartheta(\lambda_k)\bbI_{m_k}\overset{a.s.}{\to} 0,
\end{equation*}
where 
\begin{align}
	\vartheta(\lambda_k)=1+2\lambda_kc_2m_3(\lambda_k)+c_2\lambda_k^2m_3'(\lambda_k).
\end{align}
By classical CLT, we have	$[\sqrt{N}\boldsymbol{\Omega}^{\boldsymbol{F}}]_{kk}$ tends to an $m_k$-dimensional GOE (GUE) matrix under real (complex) case with the scale parameter $\lambda_k^2\vartheta$.
In fact, the scale parameter is the limit of 
\begin{align*}	&\frac{l_i^2}{N}\tr\left[\bbI+l_i\frac{1}{N}\bbY_2^*\left(\boldsymbol\Sigma_2-l_i\frac{1}{N}\bbY_2\bbY_2^*\right)^{-1}\bbY_2\right]^{2}\\
	=&\frac{l_i^2}{N}\tr\bbI_N+2\frac{l_i^3}{N}\tr \left(\boldsymbol\Sigma_2\left(\frac{1}{N}\bbY_2\bbY_2^*\right)^{-1}-l_i\bbI\right)^{-1}+\frac{l_i^4}{N}\tr \left(\boldsymbol\Sigma_2\left(\frac{1}{N}\bbY_2\bbY_2^*\right)^{-1}-l_i\bbI\right)^{-2}\\
	\overset{a.s.}{\to}&\lambda_k^2(1+2c_2\lambda_km_3(\lambda_k)+\lambda_k^2c_2m_3'(\lambda_k))=\lambda_k^2\vartheta(\lambda_k).
\end{align*}
Then we conclude that the conditional limiting distribution of $\gamma_{Nk}|\bbC_n$ equals the joint distribution of the eigenvalues of GOE (GUE) matrix with the scale parameter $1/\vartheta$.

\subsubsection{The limiting distribution of $\gamma_{Nk}$}
In this part, we will give the asymptotic distribution of $\gamma_{Nk}=(\sqrt{n}(l_i/\lambda_{Nk}-1), i\in\mathcal{J}_k)$. It is worth pointing out that the asymptotic distribution of $\gamma_{Nk}=(\sqrt{n}(l_i/\lambda_{Nk}-1), i\in\mathcal{J}_k)$ is without the condition ${\bbC_n}$.  According to (\ref{gama}), we have
\begin{align}\label{gama2}
	\gamma_{Nk}=&\sqrt{n}\frac{l_i-\lambda_{Nk}}{\lambda_{Nk}}=\sqrt{n}\frac{l_i-\psi_{\boldsymbol F}(l_i^{\bbC_n})+\psi_{\boldsymbol F}(l_i^{\bbC_n})-\psi_{\boldsymbol F}(\lambda_{nk}^{\bbC})}{\psi_{\boldsymbol F}(\lambda_{nk}^{\bbC})} \\
	=&\frac{\sqrt{n}}{\sqrt{N}}\sqrt{N}\frac{l_i-\psi_{\boldsymbol F}(l_i^{\bbC_n})}{\psi_{\boldsymbol F}(l_i^{\bbC_n})}\frac{\psi_{\boldsymbol F}(l_i^{\bbC_n})}{\psi_{\boldsymbol F}(\lambda_{nk}^{\bbC})}+\sqrt{n}\frac{l_i^{\bbC_n}-\lambda_{nk}^{\bbC}}{\lambda_{nk}^{\bbC}}\frac{\lambda_{nk}^{\bbC}}{\psi_{\boldsymbol F}(\lambda_{nk}^{\bbC})}\psi_{\boldsymbol F}'(\lambda_{nk}^{\bbC})(1+o(1)),\nonumber
\end{align} 
where the condition limiting distribution of $$\sqrt{N}\frac{l_i-\psi_{\boldsymbol F}(l_i^{\bbC_n})}{\psi_{\boldsymbol F}(l_i^{\bbC_n})}$$
is independent with the condition, so the asymptotic distribution of the first term of (\ref{gama2}) is the joint distribution of the order eigenvalues of a GOE (GUE) matrix with parameter $c_2/(c_1*\vartheta)$. According to Subsection \ref{CLT2}, it follows that the limiting distribution of $$\sqrt{n}\frac{l_i^{\bbC_n}-\lambda_{nk}^{\bbC}}{\lambda_{nk}^{\bbC}}$$ 
equals the joint distribution of the order eigenvalues of a GOE (GUE) matrix with parameter $\theta_1$ defined in (\ref{theta1}). Combining (\ref{relation3}) with (\ref{relation302}), we obtain 
$$\frac{\lambda_{nk}^{\bbC}}{\psi_{\boldsymbol F}(\lambda_{nk}^{\bbC})}=\frac{\lambda_{nk}^{\bbC}}{\lambda_{Nk}}\overset{a.s.}{\longrightarrow}1+c_2\lambda_{k}m_3(\lambda_{k}),$$
and 
\begin{align*}
	\psi_{\boldsymbol{F}}'(\lambda_{nk}^{\bbC})\overset{a.s.}{\longrightarrow}\frac{1-c_2(\lambda_{k}^{\bbC})^2m_2'(\lambda_{k}^{\bbC})}{(1+c_2\lambda_{k}m_3(\lambda_{k}))^2}.
\end{align*}
Then the asymptotic distribution of the second term of (\ref{gama2}) is the same as the eigenvalues of a GOE (GUE) matrix with parameter 
\begin{align*}
	\left[\frac{1-c_2(\lambda_{k}^{\bbC})^2m_2'(\lambda_{k}^{\bbC})}{1+c_2\lambda_{k}m_3(\lambda_{k})}\right]^2\theta_1.
\end{align*}
In summary, the limiting distribution of $\gamma_{Nk}$ is related to that of eigenvalues of the GOE (GUE) matrix with parameter 
\begin{align*}
	\frac{c_2}{c_1\cdot\vartheta}+\left[\frac{1-c_2(\lambda_{k}^{\bbC})^2m_2'(\lambda_{k}^{\bbC})}{1+c_2\lambda_{k}m_3(\lambda_{k})}\right]^2\theta_1.
\end{align*}

\subsection{Proof of Theorem \ref{cltcca}}\label{CLT4}
According to Theorem \ref{translation}, we know that there exsits a function relation between sample canonical correlation coefficients and the eigenvalues of a special noncentral Fisher matrix. The noncentral parameter matrix defined in (\ref{Noncenp}) is a random matrix, so the proof of Theorem \ref{cltcca} cannot be obtained directly be Theorem \ref{CLTnonF} and \ref{translation}. Now we will present the details of the proof.

Consider the random variable
\begin{align*}
	\gamma_k^0=\sqrt{q}\left(\frac{l_i-\Psi(\alpha_k)}{\Psi(\alpha_k)}\right),\quad \mbox{for}\ i\in\mathcal{J}_k,
\end{align*}
we have 
\begin{align*}
	\gamma_{k}^0=\sqrt{q}\left(\frac{l_i-\psi_{\bbF}\circ\psi_{\bbC}(l_i^{\boldsymbol{\Xi}})}{\psi_{\bbF}\circ\psi_{\bbC}(l_i^{\boldsymbol{\Xi}})}\frac{\psi_{\bbF}\circ\psi_{\bbC}(l_i^{\boldsymbol{\Xi}})}{\Psi(\alpha_k)}+\frac{\psi_{\bbF}\circ\psi_{\bbC}(l_i^{\boldsymbol{\Xi}})-\Psi(\alpha_k)}{\Psi(\alpha_k)}\right).
\end{align*}
Under the Assumption {d$'$} and given $\hat{\bbY}$, the limiting distribution of first term in $\gamma_{k}^0$ can be obtained by Theorem \ref{CLTnonF} and its covariance satisfies (\ref{theta}), which is independent of selection $\hat{\bbY}$. We apply the Mean value theorem to the second term in $\gamma_{k}^0$, i.e.,
\begin{align*}
	\sqrt{q}\frac{\psi_{\bbF}\circ\psi_{\bbC}(l_i^{\boldsymbol{\Xi}})\!-\!\Psi(\alpha_k)}{\Psi(\alpha_k)}=\sqrt{n}\frac{l_i^{\boldsymbol{\Xi}}\!-\!\psi_{\boldsymbol{\Xi}}(f(\alpha_k))}{\psi_{\boldsymbol{\Xi}}(f(\alpha_k))}\frac{\sqrt{q}}{\sqrt{n}}\frac{\psi_{\boldsymbol{\Xi}}(f(\alpha_k))}{\Psi(\alpha_k)}\psi_{\bbF}'(\xi_1)\psi_{\bbC}'(\xi_2),
\end{align*}
where $\xi_1 \in (\Psi_{\bbC}(\alpha_k), \Psi_{\bbC}(l_i^{\boldsymbol{\Xi}}))$ or $(\Psi_{\bbC}(l_i^{\boldsymbol{\Xi}}),\Psi_{\bbC}(\alpha_k))$, and $\xi_2\in (\psi_{\boldsymbol{\Xi}}(f(\alpha_k)),l_i^{\boldsymbol{\Xi}})$ or $(l_i^{\boldsymbol{\Xi}}, \psi_{\boldsymbol{\Xi}}(f(\alpha_k)))$. By Theorem \ref{limittuS1} and \ref{limitsF}, we have 
\begin{align*}
	&\psi_F'(\xi_1)-\psi_{\bbF}'(\Psi_{\bbC}(\alpha_k))\overset{a.s.}{\to}0\\
	&\psi_C'(\xi_2)-\psi_{\bbC}'(\psi_{\boldsymbol{\Xi}}(f(\alpha_k)))\overset{a.s.}{\to}0
\end{align*}
where
\begin{align*}
&\psi_{\bbF}'(\Psi_{\bbC}(\alpha_k))\overset{\triangle}{=}\frac{1-c_4\Psi_{\bbC}^2(\alpha_k)m_{\bbC}'(\Psi_{\bbC}(\alpha_k))}{[1+c_4\Psi_{\bbC}(\alpha_k)m_{\bbC}(\Psi_{\bbC}(\alpha_k))]^2},\\
&\psi_{\bbC}'(\psi_{\boldsymbol{\Xi}}(f(\alpha_k)))\overset{\triangle}{=}\left(1-c_3\int\frac{dF_{mp}^{p/n,H}(t)}{t-\psi_{\boldsymbol{\Xi}}(f(\alpha_k))}\right)^2\\
&\quad\quad-2\psi_{\boldsymbol{\Xi}}(f(\alpha_k))\left(1-c_3\int\frac{dF_{mp}^{p/n,H}(t)}{t-\psi_{\boldsymbol{\Xi}}(f(\alpha_k))}\right)c_3\int\frac{dF_{mp}^{p/n,H}(t)}{(t-\psi_{\boldsymbol{\Xi}}(f(\alpha_k)))^2}\\
&\quad\quad-(1-c_3)c_3\int\frac{dF_{mp}^{p/n,H}(t)}{(t-\psi_{\boldsymbol{\Xi}}(f(\alpha_k)))^2}.
\end{align*}
And 
\begin{align*}
	\sqrt{n}\frac{l_i^{\boldsymbol{\Xi}}\!-\!\psi_{\boldsymbol{\Xi}}(f(\alpha_k))}{\psi_{\boldsymbol{\Xi}}(f(\alpha_k))}\times \sqrt{\frac{2}{(\psi_{\boldsymbol{\Xi}}^2(f(\alpha_k))\underline{m}'(\psi_{\boldsymbol{\Xi}}(f(\alpha_k)))}}\overset{d}{\to}N(0,1),
\end{align*}
where
\begin{align*}
	\underline{m}(\psi_{\boldsymbol{\Xi}}(f(\alpha_k)))=-\frac{1-p/n}{\psi_{\boldsymbol{\Xi}}(f(\alpha_k))}+p/n\int\frac{dF_{mp}^{p/n,H}(t)}{t-\psi_{\boldsymbol{\Xi}}(f(\alpha_k))}
\end{align*}
Then the covariance function of $\gamma_k^0$ satisfies 
\begin{align}\label{Cov1}
	\eta_3\!+\!\frac{c_2}{c_1\cdot\eta_2}\!+\!\left[\frac{1-c_2(\lambda_{k}^{\bbC})^2m_2'(\lambda_{k}^{\bbC})}{1\!+\!c_2\lambda_{k}m_3(\lambda_{k})}\right]^2\eta_1.
\end{align}
where
\begin{align*}
	\eta_3=\frac{q}{n}\frac{(\psi_{\bbF}'(\Psi_{\bbC}(\alpha_k)))^2(\psi_{\bbC}'(\psi_{\boldsymbol{\Xi}}(f(\alpha_k))))^2}{\psi_{\boldsymbol{\Xi}}^2(f(\alpha_k))\underline{m}'(\psi_{\boldsymbol{\Xi}}(f(\alpha_k)))}\frac{\psi_{\boldsymbol{\Xi}}^2(f(\alpha_k))}{\Psi^2(\alpha_k)}
\end{align*}

According to Delta method, we rewrite
\begin{align*}
	\gamma_k=\sqrt{q}\frac{\lambda_i^2-t(\alpha_k)}{t(\alpha_k)}=\sqrt{q}\frac{g^{-1}(l_i)-g^{-1}(\Psi(\alpha_k))}{\Psi(\alpha_k)}\frac{\Psi(\alpha_k)}{t(\alpha_k)}
\end{align*} 
then the covariance function of $\gamma_k$ will be 
\begin{align*}
	(\ref{Cov1})*\left(\frac{c_4}{[1+c_4\Psi(\alpha_k)]^2}\right)^2*\frac{\Psi^2(\alpha_k)}{t^2(\alpha_k)}
\end{align*}

\end{document}